\documentclass[reqno,11pt]{amsart}
\usepackage{amsfonts}

\usepackage{graphicx}
\usepackage{pstricks}
\usepackage{amsmath}
\usepackage{amsmath}
\usepackage{amsxtra}

\theoremstyle{plain}
\newtheorem{theorem}{Theorem}[section]
\newtheorem{lemma}[theorem]{Lemma}

\newtheorem{proposition}[theorem]{Proposition}
\newtheorem{corollary}[theorem]{Corollary}
\newtheorem{conjecture}[theorem]{Conjecture}
\newtheorem{definition}[theorem]{Definition}

\newtheorem{problem}[theorem]{Problem}

\theoremstyle{definition}
\newtheorem{notation}[theorem]{Notation}
\newtheorem{example}[theorem]{Example}
\newtheorem{remark}[theorem]{Remark}

\numberwithin{equation}{section}
\begin{document}
\title{Jointly hyponormal pairs of commuting subnormal operators need not be
jointly subnormal}
\author{Ra\'{u}l E. Curto}
\address{Department of Mathematics, The University of Iowa, Iowa City, Iowa
52242}
\email{rcurto@math.uiowa.edu}
\urladdr{http://www.math.uiowa.edu/\symbol{126}rcurto/}
\author{Jasang Yoon}
\address{Department of Mathematics, Iowa State University, Ames, Iowa 50011}
\email{yoon@math.uiowa.edu}
\urladdr{http://math.uiowa.edu/\symbol{126}yoon/}
\thanks{Research partially supported by NSF Grant DMS-0099357. \ }

\begin{abstract}
We construct three different families of commuting pairs of subnormal
operators, jointly hyponormal but not admitting commuting normal extensions.
\ Each such family can be used to answer in the negative a 1988 conjecture
of R. Curto, P. Muhly and J. Xia. \ We also obtain a sufficient condition
under which joint hyponormality does imply joint subnormality.
\end{abstract}

\maketitle

\section{\label{Int}Introduction}

Let $\mathcal{H}$ be a complex Hilbert space and let $\mathcal{B}(\mathcal{H}%
)$ denote the algebra of bounded linear operators on $\mathcal{H}$. $\ $For $%
S,T\in \mathcal{B}(\mathcal{H})$ let $[S,T]:=ST-TS$. \ We say that an $n$%
-tuple $\mathbf{T}=(T_{1},\ldots ,T_{n})$ of operators on $\mathcal{H}$ is
(jointly) \textit{hyponormal} if the operator matrix 
\begin{equation*}
\lbrack \mathbf{T}^{\ast },\mathbf{T]:=}\left( 
\begin{array}{llll}
\lbrack T_{1}^{\ast },T_{1}] & [T_{2}^{\ast },T_{1}] & \cdots & [T_{n}^{\ast
},T_{1}] \\ 
\lbrack T_{1}^{\ast },T_{2}] & [T_{2}^{\ast },T_{2}] & \cdots & [T_{n}^{\ast
},T_{2}] \\ 
\text{ \thinspace \thinspace \quad }\vdots & \text{ \thinspace \thinspace
\quad }\vdots & \cdots & \text{ \thinspace \thinspace \quad }\vdots \\ 
\lbrack T_{1}^{\ast },T_{n}] & [T_{2}^{\ast },T_{n}] & \cdots & [T_{n}^{\ast
},T_{n}]%
\end{array}%
\right)
\end{equation*}%
is positive on the direct sum of $n$ copies of $\mathcal{H}$ (cf. \cite{Ath}%
, \cite{CMX}). \ The $n$-tuple $\mathbf{T}$ is said to be \textit{normal} if 
$\mathbf{T}$ is commuting and each $T_{i}$ is normal, and $\mathbf{T}$ is 
\textit{subnormal }if $\mathbf{T}$ is the restriction of a normal $n$-tuple
to a common invariant subspace. \ Clearly, normal $\Rightarrow $ subnormal $%
\Rightarrow $ hyponormal. \ The Bram-Halmos criterion states that an
operator $T\in \mathcal{B}(\mathcal{H})$ is subnormal if and only the $k$%
-tuple $(T,T^{2},...,T^{k})$ is hyponormal for all $k\geq 1$.

For $\alpha \equiv \{\alpha _{n}\}_{n=0}^{\infty }$ a bounded sequence of
positive real numbers (called \textit{weights}), let $W_{\alpha }:\ell ^{2}(%
\mathbb{Z}_{+})\rightarrow \ell ^{2}(\mathbb{Z}_{+})$ be the associated
unilateral weighted shift, defined by $W_{\alpha }e_{n}:=\alpha
_{n}e_{n+1}\;($all $n\geq 0)$, where $\{e_{n}\}_{n=0}^{\infty }$ is the
canonical orthonormal basis in $\ell ^{2}(\mathbb{Z}_{+}).$ \ The moments of 
$\alpha $ are given as 
\begin{equation*}
\gamma _{k}\equiv \gamma _{k}(\alpha ):=\left\{ 
\begin{array}{cc}
1 & \text{if }k=0 \\ 
\alpha _{0}^{2}\cdot ...\cdot \alpha _{k-1}^{2} & \text{if }k>0%
\end{array}%
\right\} .
\end{equation*}%
It is easy to see that $W_{\alpha }$ is never normal, and that it is
hyponormal if and only if $\alpha _{0}\leq \alpha _{1}\leq ...$. \
Similarly, consider double-indexed positive bounded sequences $\alpha _{%
\mathbf{k}},\beta _{\mathbf{k}}\in \ell ^{\infty }(\mathbb{Z}_{+}^{2})$ , $%
\mathbf{k}\equiv (k_{1},k_{2})\in \mathbb{Z}_{+}^{2}:=\mathbb{Z}_{+}\times 
\mathbb{Z}_{+}$ and let $\ell ^{2}(\mathbb{Z}_{+}^{2})$\ be the Hilbert
space of square-summable complex sequences indexed by $\mathbb{Z}_{+}^{2}$.
\ (Recall that $\ell ^{2}(\mathbb{Z}_{+}^{2})$ is canonically isometrically
isomorphic to $\ell ^{2}(\mathbb{Z}_{+})\bigotimes \ell ^{2}(\mathbb{Z}_{+})$%
.) \ We define the $2$-variable weighted shift $\mathbf{T}\equiv
(T_{1},T_{2})$\ by 
\begin{equation*}
T_{1}e_{\mathbf{k}}:=\alpha _{\mathbf{k}}e_{\mathbf{k+}\varepsilon _{1}}
\end{equation*}%
\begin{equation*}
T_{2}e_{\mathbf{k}}:=\beta _{\mathbf{k}}e_{\mathbf{k+}\varepsilon _{2}},
\end{equation*}%
where $\mathbf{\varepsilon }_{1}:=(1,0)$ and $\mathbf{\varepsilon }%
_{2}:=(0,1)$. \ Clearly, 
\begin{equation}
T_{1}T_{2}=T_{2}T_{1}\Longleftrightarrow \beta _{\mathbf{k+}\varepsilon
_{1}}\alpha _{\mathbf{k}}=\alpha _{\mathbf{k+}\varepsilon _{2}}\beta _{%
\mathbf{k}}\;\;(\text{all }\mathbf{k}).  \label{commuting}
\end{equation}%
In an entirely similar way one can define multivariable weighted shifts. \
Trivially, a pair of unilateral weighted shifts $W_{a}$ and $W_{\beta }$
gives rise to a $2$-variable weighted shift $\mathbf{T}\equiv (T_{1},T_{2})$%
, if we let $\alpha _{(k_{1},k_{2})}:=\alpha _{k_{1}}$ and $\beta
_{(k_{1},k_{2})}:=\beta _{k_{2}}\;$(all $k_{1},k_{2}\in \mathbb{Z}_{+}^{2}$%
). \ In this case, $\mathbf{T}$ is subnormal (resp. hyponormal) if and only
if so are $T_{1}$ and $T_{2}$; in fact, under the canonical identification
of $\ell ^{2}(\mathbb{Z}_{+}^{2})$ and $\ell ^{2}(\mathbb{Z}_{+})\bigotimes
\ell ^{2}(\mathbb{Z}_{+})$, $T_{1}\cong W_{a}\bigotimes I$ and $T_{2}\cong
I\bigotimes W_{\beta }$, and $\mathbf{T}$ is also doubly commuting. \ For
this reason, we do not focus attention on shifts of this type, and use them
only when the above mentioned triviality is desirable or needed. \ 

We now recall a well known characterization of subnormality for single
variable weighted shifts, due to C. Berger (cf. \cite[III.8.16]{Con}): \ $%
W_{\alpha }$ is subnormal if and only if there exists a probability measure $%
\xi $ supported in $[0,\left\| W_{\alpha }\right\| ^{2}]$ such that $\gamma
_{k}(\alpha ):=\alpha _{0}^{2}\cdot ...\cdot \alpha _{k-1}^{2}=\int
t^{k}\;d\xi (t)\;\;(k\geq 1)$. \ If $W_{\alpha }$ is subnormal, and if for $%
h\geq 1$ we let $\mathcal{M}_{h}:=\bigvee \{e_{n}:n\geq h\}$ denote the
invariant subspace obtained by removing the first $h$ vectors in the
canonical orthonormal basis of $\ell ^{2}(\mathbb{Z}_{+})$, then the Berger
measure of $W_{\alpha }|_{\mathcal{M}_{h}}$ is $\frac{1}{\gamma _{h}}%
t^{h}d\xi (t)$.

An important class of subnormal weighted shifts is obtained by considering
measures $\mu $ with exactly two atoms $t_{0}$ and $t_{1}$. \ These shifts
arise naturally in the Subnormal Completion Problem \cite{RGWSII} and in the
theory of truncated moment problems (cf. \cite{Houston}, \cite{tcmp1}). \
For $t_{0},t_{1}\in \mathbb{R}_{+}$, $t_{0}<t_{1}$, and $\rho _{0},\rho
_{1}>0$, the moments of the $2$-atomic measure $\xi :=\rho _{0}\delta
_{t_{0}}+\rho _{1}\delta _{t_{1}}$ (here $\delta _{p}$ denotes the
point-mass probability measure with support the singleton $\{p\}$) satisfy
the $2$-step recursive relation $\gamma _{n+2}=\varphi _{0}\gamma
_{n}+\varphi _{1}\gamma _{n+1}\;(n\geq 0)$; at the weight level, this can be
written as $\alpha _{n+1}^{2}=\frac{\varphi _{0}}{\alpha _{n}^{2}}+\varphi
_{1}\;(n\geq 0)$. \ More generally, any finitely atomic Berger measure
corresponds to a recursively generated weighted shift (i.e., one whose
moments satisfy an $r$-step recursive relation); in fact, $r=$ card supp $%
\xi $. \ In the special case of $r=2$, the theory of recursively generated
weighted shifts makes contact with the work of J. Stampfli in \cite{Sta}, in
which he proved that given three positive numbers $\alpha _{0}<\alpha
_{1}<\alpha _{2}$, it is always possible to find a subnormal weighted shift,
denoted $W_{(\alpha _{0},\alpha _{1},\alpha _{2})\symbol{94}}$, whose first
three weights are $\alpha _{0},\alpha _{1}$ and $\alpha _{2}$. \ In this
case, the coefficients of recursion (cf. \cite[Example 3.12]{RGWSI}, %
\cite[Section 3]{RGWSII}, \cite[Section 1, p. 81]{OTAMP}) are given by 
\begin{equation}
\varphi _{0}=-\frac{\alpha _{0}^{2}\alpha _{1}^{2}(\alpha _{2}^{2}-\alpha
_{1}^{2})}{\alpha _{1}^{2}-\alpha _{0}^{2}}\text{ and }\varphi _{1}=\frac{%
\alpha _{1}^{2}(\alpha _{2}^{2}-\alpha _{0}^{2})}{\alpha _{1}^{2}-\alpha
_{0}^{2}}\text{,}  \label{phieq}
\end{equation}%
the atoms $t_{0}$ and $t_{1}$ are the roots of the equation 
\begin{equation}
t^{2}-(\varphi _{0}+\varphi _{1}t)=0,  \label{teq}
\end{equation}%
and the densities $\rho _{0}$ and $\rho _{1}$ uniquely solve the $2\times 2$
system of equations 
\begin{equation}
\left\{ 
\begin{array}{ccc}
\rho _{0}+\rho _{1} & = & 1 \\ 
\rho _{0}t_{0}+\rho _{1}t_{1} & = & \alpha _{0}^{2}%
\end{array}%
\right. .  \label{rhoeq}
\end{equation}

We also recall the notion of moment of order $\mathbf{k}$ for a pair $%
(\alpha ,\beta )$ satisfying (\ref{commuting}). \ Given $\mathbf{k}\in 
\mathbb{Z}_{+}^{2}$, the moment of $(\alpha ,\beta )$ of order $\mathbf{k}$
is 
\begin{equation*}
\gamma _{\mathbf{k}}\equiv \gamma _{\mathbf{k}}(\alpha ,\beta ):=\left\{ 
\begin{array}{cc}
1 & \text{if }\mathbf{k}=0 \\ 
\alpha _{(0,0)}^{2}\cdot ...\cdot \alpha _{(k_{1}-1,0)}^{2}\cdot \beta
_{(k_{1},0)}^{2}\cdot ...\cdot \beta _{(k_{1},k_{2}-1)}^{2} & \text{if }%
\mathbf{k\in }\mathbb{Z}_{+}^{2}\text{, }\mathbf{k}\neq 0%
\end{array}%
\right\} .
\end{equation*}%
We remark that, due to the commutativity condition (\ref{commuting}), $%
\gamma _{\mathbf{k}}$ can be computed using any nondecreasing path from $%
(0,0)$ to $(k_{1},k_{2})$.

\begin{theorem}
\label{Berger}(Berger's Theorem, $2$-variable case) (\cite{JeLu}) $\ $A $2$%
-variable weighted shift $\mathbf{T}\equiv (T_{1},T_{2})$ admits a commuting
normal extension if and only if there is a probability measure $\mu $
defined on the $2$-dimensional rectangle $R=[0,a_{1}]\times \lbrack 0,a_{2}]$
($a_{i}:=\left\| T_{i}\right\| ^{2}$) such that $\gamma _{\mathbf{k}%
}=\iint_{R}\mathbf{t}^{\mathbf{k}}d\mu (\mathbf{t}):=%
\iint_{R}t_{1}^{k_{1}}t_{2}^{k_{2}}\;d\mu (t_{1},t_{2})$ \ $($all $\mathbf{%
k\in }\mathbb{Z}_{+}^{2}$).
\end{theorem}

Clearly, each component $T_{i}$ of a subnormal $2$-variable weighted shift $%
\mathbf{T}\equiv (T_{1},T_{2})$ must be subnormal. \ For instance, $%
T_{1}\cong \bigoplus_{j=0}^{\infty }W_{\alpha ^{(j)}}$, where $\alpha
_{i}^{(j)}:=\alpha _{(i,j)}$, so that $W_{\alpha ^{(j)}}$ has associated
Berger measure $d\nu _{j}(t_{1}):=\frac{1}{\gamma _{(0,j)}}%
\int_{[0,a_{2}]}t_{2}^{j}d\Phi _{t_{1}}(t_{2})$, where $d\mu
(t_{1},t_{2})\equiv d\Phi _{t_{1}}(t_{2})d\eta (t_{1})$ is the canonical
disintegration of $\mu $ by horizontal slices. \ On the other hand, if we
only know that each of $T_{1}$, $T_{2}$ is subnormal, and that they commute,
the following problem is natural.

\begin{problem}
(Lifting Problem for Commuting Subnormals) \ Find necessary and sufficient
conditions on $T_{1}$ and $T_{2}$ to guarantee the subnormality of $\mathbf{T%
}\equiv (T_{1},T_{2})$. \ 
\end{problem}

It is well known that the above mentioned necessary conditions do not
suffice (cf.\cite{bridge}). \ In terms of the \textit{marginal} measures,
the problem can be phrased as a reconstruction-of-measure problem, that is,
under what conditions on the single variable measures $\{\nu
_{j}\}_{j=0}^{\infty }$ and $\{\omega _{i}\}_{i=0}^{\infty }$ associated
with $T_{1}$ and $T_{2}$, respectively, does there exist a $2$-variable
measure $\mu $ correctly interpolating all the powers $%
t_{1}^{k_{1}}t_{2}^{k_{2}}\;(k_{1},k_{2}\geq 0)$. $\ $

To detect hyponormality for $2$-variable weighted shifts, there is a simple
criterion involving a base point $\mathbf{k}$ in $\mathbb{Z}_{+}^{2}$ and
its five neighboring points in $\mathbf{k}+\mathbb{Z}_{+}^{2}$ at path
distance at most $2$ (cf. Figure \ref{Figure 0}).

\begin{theorem}
(\cite{bridge})\label{joint hypo} \ (Six-point Test) \ Let $\mathbf{T\equiv (%
}T_{1},T_{2})$ be a $2$-variable weighted shift, with weight sequences $%
\alpha $ and $\beta $. \ Then 
\begin{equation*}
\lbrack \mathbf{T}^{\ast },\mathbf{T]\geq }0\Leftrightarrow (([T_{j}^{\ast
},T_{i}]e_{\mathbf{k+}\varepsilon _{j}},e_{\mathbf{k+}\varepsilon
_{i}}))_{i,j=1}^{2}\geq 0\text{ (all }\mathbf{k}\in \mathbf{Z}_{+}^{2}\text{)%
}
\end{equation*}%
\begin{equation*}
\Leftrightarrow \left( 
\begin{array}{cc}
\alpha _{\mathbf{k}+\mathbf{\varepsilon }_{1}}^{2}-\alpha _{\mathbf{k}}^{2}
& \alpha _{\mathbf{k}+\mathbf{\varepsilon }_{2}}\beta _{\mathbf{k}+\mathbf{%
\varepsilon }_{1}}-\alpha _{\mathbf{k}}\beta _{\mathbf{k}} \\ 
\alpha _{\mathbf{k}+\mathbf{\varepsilon }_{2}}\beta _{\mathbf{k}+\mathbf{%
\varepsilon }_{1}}-\alpha _{\mathbf{k}}\beta _{\mathbf{k}} & \beta _{\mathbf{%
k}+\mathbf{\varepsilon }_{2}}^{2}-\beta _{\mathbf{k}}^{2}%
\end{array}%
\right) \geq 0\text{ (all }\mathbf{k}\in \mathbf{Z}_{+}^{2}\text{).}
\end{equation*}
\end{theorem}

\setlength{\unitlength}{1mm} \psset{unit=1mm}

\begin{figure}[th]
\begin{center}
\begin{picture}(90,60)

\psline(20,10)(60,10)
\psline(20,30)(40,30)
\psline(20,10)(20,50)
\psline(40,10)(40,30)

\put(20,10){\pscircle*(0,0){1}}
\put(40,10){\pscircle*(0,0){1}}
\put(60,10){\pscircle*(0,0){1}}
\put(20,30){\pscircle*(0,0){1}}
\put(40,30){\pscircle*(0,0){1}}
\put(20,50){\pscircle*(0,0){1}}

\put(13,6){\footnotesize{$(k_1,k_2)$}}
\put(31,6){\footnotesize{$(k_1+1,k_2)$}}
\put(54,6){\footnotesize{$(k_1+2,k_2)$}}

\put(27,12){\footnotesize{$\alpha_{k_1,k_2}$}}
\put(47,12){\footnotesize{$\alpha_{k_1+1,k_2}$}}

\put(27,32){\footnotesize{$\alpha_{k_1,k_2+1}$}}

\put(3,29){\footnotesize{$(k_1,k_2+1)$}}
\put(3,49){\footnotesize{$(k_1, k_2+2)$}}
\put(41,29){\footnotesize{$(k_1+1, k_2+1)$}}

\put(20,19){\footnotesize{$\beta_{k_1,k_2}$}}
\put(20,39){\footnotesize{$\beta_{k_1,k_2+1}$}}

\put(40,19){\footnotesize{$\beta_{k_1+1,k_2}$}}

\end{picture}
\end{center}
\caption{{}Weight diagram used in the Six-point Test}
\label{Figure 0}
\end{figure}
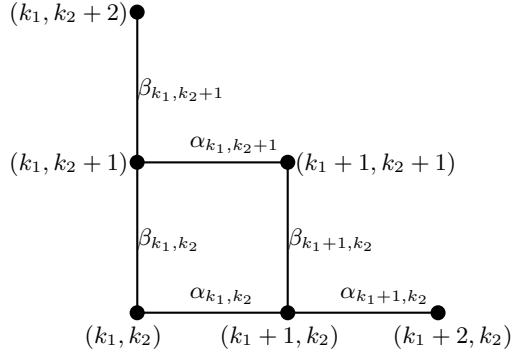

Unlike the single variable case, in which there is a clear separation
between hyponormality and subnormality (cf. \cite{RGWSII}, \cite{OTAMP}, %
\cite{CuLe3}), much less is known about the multivariable case. \ In this
paper we will construct three conceptually different families of
counterexamples to the following conjecture.

\begin{conjecture}
\label{conjecture}(\cite{CMX}) \ Let $\mathbf{T\equiv(}T_{1},T_{2})$ be a
pair of commuting subnormal operators on $\mathcal{H}$. $\ $Then $\mathbf{T}$
is subnormal if and only if $\mathbf{T}$ is hyponormal.
\end{conjecture}

We mention that M. Dritschel and S. McCullough, working independently, have
been able to obtain a separate example (\cite{DrMcC}). \ We shall see in
Section \ref{fourth} that their example is a special case of a general
construction that produces nonsubnormal hyponormal pairs with $T_{1}\cong
T_{2}$.

We now formulate an improved version of a result due to the first-named
author.

\begin{proposition}
\label{backward}(Subnormal backward extension of a $1$-variable weighted
shift) (cf \cite{QHWS}) \ Let $T$ be a weighted shift whose restriction $T_{%
\mathcal{M}}$ to $\mathcal{M}:=\vee \{e_{1},e_{2},\cdots \}$ is subnormal,
with associated measure $\mu _{\mathcal{M}}.$ \ Then $T$ is subnormal (with
associated measure $\mu $) if and only if\newline
(i) $\ \frac{1}{t}\in L^{1}(\mu _{\mathcal{M}})$\newline
(ii) $\ \alpha _{0}^{2}\leq (\left\| \frac{1}{t}\right\| _{L^{1}(\mu _{%
\mathcal{M}})})^{-1}$\newline
In this case, $d\mu (t)=\frac{\alpha _{0}^{2}}{t}d\mu _{\mathcal{M}%
}(t)+(1-\alpha _{0}^{2}\left\| \frac{1}{t}\right\| _{L^{1}(\mu _{\mathcal{M}%
})})d\delta _{0}(t)$, where $\delta _{0}$ denotes Dirac measure at $0$. \ In
particular, $T$ is never subnormal when $\mu _{\mathcal{M}}(\{0\})>0$. \ 
\end{proposition}

\begin{proof}
$\Rightarrow )$ \ We first observe that the moments of $T$ and $T_{\mathcal{M%
}}$ are related by the equation 
\begin{equation*}
\gamma _{k}(T_{\mathcal{M}})\equiv \alpha _{1}^{2}\cdots \alpha _{k}^{2}=%
\frac{\gamma _{k+1}\left( T\right) }{\alpha _{0}^{2}}
\end{equation*}
so that 
\begin{equation*}
\frac{1}{\alpha _{0}^{2}}\int t^{k+1}d\mu (t)=\int t^{k}d\mu _{\mathcal{M}%
}(t)\;\;(\text{all }k\geq 0),
\end{equation*}
that is, $td\mu (t)=\alpha _{0}^{2}d\mu _{\mathcal{M}}(t)$. \ \ It follows
at once that 
\begin{equation*}
d\mu (t)=\lambda d\delta _{0}(t)+\frac{\alpha _{0}^{2}}{t}d\mu _{\mathcal{M}%
}(s),
\end{equation*}
where $\lambda \geq 0.$ \ Since $\int d\mu =1,$ we must have $\frac{1}{t}\in
L^{1}(\mu _{\mathcal{M}})$ and $\alpha _{0}^{2}\left\| \frac{1}{t}\right\|
_{L^{1}(\mu _{\mathcal{M}})}\leq 1$. \ Finally, it is straightforward to
verify that $\lambda =(1-\alpha _{0}^{2}\left\| \frac{1}{t}\right\|
_{L^{1}(\mu _{\mathcal{M}})})$.

$\Leftarrow )$ \ Let%
\begin{equation*}
d\mu (t):=\frac{\alpha _{0}^{2}}{t}d\mu _{\mathcal{M}}(t)+(1-\alpha
_{0}^{2}\left\| \frac{1}{t}\right\| _{L^{1}(\mu _{\mathcal{M}})})d\delta
_{0}(t).
\end{equation*}%
By hypotheses, $\mu $ is a positive Borel measure on $[0,\left\| T\right\|
^{2}]$. \ Moreover, 
\begin{equation*}
\int d\mu =\alpha _{0}^{2}\int \frac{1}{t}d\mu _{\mathcal{M}}+(1-\alpha
_{0}^{2}\left\| \frac{1}{t}\right\| _{L^{1}(\mu _{\mathcal{M}})})\int
d\delta _{0}=1,
\end{equation*}%
and for $k\geq 1$,%
\begin{eqnarray*}
\int t^{k}d\mu (t) &=&\alpha _{0}^{2}\int t^{k}\frac{1}{t}d\mu _{\mathcal{M}%
}(t)+(1-\alpha _{0}^{2}\left\| \frac{1}{t}\right\| _{L^{1}(\mu _{\mathcal{M}%
})})\int t^{k}d\delta _{0}(t) \\
&=&\alpha _{0}^{2}\int t^{k-1}d\mu _{\mathcal{M}}(t)=\alpha _{0}^{2}\gamma
_{k-1}(T_{\mathcal{M}})=\gamma _{k}(T)\text{.}
\end{eqnarray*}
Therefore, $T$ is subnormal, with Berger measure $\mu $.
\end{proof}

\begin{notation}
The maximum possible value for $\alpha _{0}$ in Proposition \ref{backward},
namely $(\left\| \frac{1}{t}\right\| _{L^{1}(\mu _{\mathcal{M}})})^{-1}$,
will be denoted by $\alpha _{ext}\equiv \alpha _{ext}(\mu _{\mathcal{M}})$.
\ More generally, given a ($1$-variable) subnormal weighted shift $W_{\eta }$
with weight sequence $\eta _{1}\leq \eta _{2}\leq ...$ and Berger measure $%
\nu $, we let 
\begin{equation*}
\eta _{ext}:=\left\{ 
\begin{array}{cc}
0 & \text{if }\frac{1}{t}\notin L^{1}(\nu ) \\ 
(\left\| \frac{1}{t}\right\| _{L^{1}(\nu )})^{-1} & \text{if }\frac{1}{t}\in
L^{1}(\nu )%
\end{array}%
\right. .
\end{equation*}%
Observe that when the weight sequence $\eta $ is strictly increasing and $%
\frac{1}{t}\in L^{1}(\nu )$, we must necessarily have $\eta _{ext}<\eta _{1}$%
, by \cite[Theorem 6]{Sta}. \ On occasion, we will write $shift(\alpha
_{0},\alpha _{1},...)$ to denote the weighted shift with weight sequence $%
\{\alpha _{k}\}_{k=0}^{\infty }$. \ We also denote by $U_{+}:=shift(1,1,...)$
the (unweighted) unilateral shift, and for $0<a<1$ we let $%
S_{a}:=shift\{a,1,1,...)$. \ Observe that the Berger measures of $U_{+}$ and 
$S_{a}$ are $\delta _{1}$ and $(1-a^{2})\delta _{0}+a^{2}\delta _{1}$,
respectively, where $\delta _{p}$ denotes the point-mass probability measure
with support the singleton $\{p\}$. \ Finally, we let $B_{+}$ denote the
Bergman shift, whose Berger measure is Lebesgue measure on the interval $%
[0,1]$; the weights of $B_{+}$ are given by the formula $\alpha _{n}:=\sqrt{%
\frac{n+1}{n+2}}\;(n\geq 0)$.
\end{notation}

We conclude this section with a result that will be needed in Section \ref%
{secondcounter}.

\begin{lemma}
\label{abclem}(cf. \cite[Theorem 3.10]{RGWSII}) \ For $0<\alpha _{0}<\alpha
_{1}<\alpha _{2}$, let $W_{(\alpha _{0},\alpha _{1},\alpha _{2})\symbol{94}}$
be the weighted shift described by (\ref{phieq}), (\ref{teq}) and (\ref%
{rhoeq}). \ Consider now $W_{\eta }:=shift(\alpha _{1},\alpha _{2},...)$,
that is, $W_{\eta }$ is the restriction of $W_{(\alpha _{0},\alpha
_{1},\alpha _{2})\symbol{94}}$ to $\mathcal{M}$. \ Then $\eta _{ext}=\alpha
_{0}$.
\end{lemma}

\textit{Acknowledgment}. \ All the examples, and the basic construction in
Section \ref{secondcounter} were obtained using calculations with the
software tool \textit{Mathematica \cite{Wol}.}

\section{\label{secondcounter}The First Family of Counterexamples}

\textbf{Construction of the family}. \ Let $0<a,b<1$ and let $\{\xi
_{k}\}_{k=0}^{\infty }$ and $\{\eta _{k}\}_{k=0}^{\infty }$ be two strictly
increasing weight sequences. \ Consider the $2$-variable weighted shift $%
\mathbf{T\equiv (}T_{1},T_{2})$ on $\ell ^{2}(\mathbb{Z}_{+}^{2})$ given by
the double-indexed weight sequences 
\begin{equation}
\alpha (\mathbf{k}):=\left\{ 
\begin{tabular}{ll}
$\xi _{k_{1}}$ & $\text{if }k_{1}\geq 1$ or $k_{2}\geq 1$ \\ 
$a$ & $\text{if }k_{1}=0\text{ and }k_{2}=0$%
\end{tabular}%
\right.  \label{eqnew1}
\end{equation}%
and 
\begin{equation}
\beta (\mathbf{k}):=\left\{ 
\begin{tabular}{ll}
$\eta _{k_{2}}$ & $\text{if }k_{1}\geq 1$ or $k_{2}\geq 1$ \\ 
$b$ & $\text{if }k_{1}=0\text{ and }k_{2}=0$%
\end{tabular}%
\right. .  \label{eqnew2}
\end{equation}%
where $W_{\xi }$ and $W_{\eta }$ are two single-variable subnormal weighted
shifts with Berger measures $\nu $ and $\omega $, resp., and 
\begin{equation}
a\eta _{0}=b\xi _{0}  \label{commuting2}
\end{equation}%
(to guarantee the commutativity of $T_{1}$ and $T_{2}$ (\ref{commuting})). \ 
$\mathbf{T}$ can be represented by the following weight diagram (Figure \ref%
{Figure 1}). \ It is then clear that $T_{1}$ and $T_{2}$ are subnormal
provided $a\leq \xi _{ext}(\nu _{\mathcal{M}})$ and $b\leq \eta
_{ext}(\omega _{\mathcal{M}})$, where, as usual, $\mathcal{M}:=\vee
\{e_{1},e_{2,}\cdots \}$; in particular, $a<\xi _{1}$ and $b<\eta _{1}$.

\setlength{\unitlength}{1mm} \psset{unit=1mm}

\begin{figure}[th]
\begin{center}
\begin{picture}(140,138)

\psline{->}(20,20)(135,20)
\psline(20,40)(125,40)
\psline(20,60)(125,60)
\psline(20,80)(125,80)
\psline(20,100)(125,100)
\psline(20,120)(125,120)
\psline{->}(20,20)(20,135)
\psline(40,20)(40,125)
\psline(60,20)(60,125)
\psline(80,20)(80,125)
\psline(100,20)(100,125)
\psline(120,20)(120,125)

\put(11,16){\footnotesize{$(0,0)$}}
\put(35,16){\footnotesize{$(1,0)$}}
\put(55,16){\footnotesize{$(2,0)$}}
\put(78,16){\footnotesize{$\cdots$}}
\put(95,16){\footnotesize{$(n,0)$}}
\put(115,16){\footnotesize{$(n+1,0)$}}

\put(27,21){\footnotesize{$a$}}
\put(47,21){\footnotesize{$\xi_{1}$}}
\put(67,21){\footnotesize{$\xi_{2}$}}
\put(87,21){\footnotesize{$\cdots$}}
\put(107,21){\footnotesize{$\xi_{n}$}}
\put(124,21){\footnotesize{$\cdots$}}

\put(27,41){\footnotesize{$\xi_{0}$}}
\put(47,41){\footnotesize{$\xi_{1}$}}
\put(67,41){\footnotesize{$\xi_{2}$}}
\put(87,41){\footnotesize{$\cdots$}}
\put(107,41){\footnotesize{$\xi_{n}$}}
\put(124,41){\footnotesize{$\cdots$}}

\put(27,61){\footnotesize{$\xi_{0}$}}
\put(47,61){\footnotesize{$\xi_{1}$}}
\put(67,61){\footnotesize{$\xi_{2}$}}
\put(87,61){\footnotesize{$\cdots$}}
\put(107,61){\footnotesize{$\xi_{n}$}}
\put(124,61){\footnotesize{$\cdots$}}

\put(27,81){\footnotesize{$\cdots$}}
\put(47,81){\footnotesize{$\cdots$}}
\put(67,81){\footnotesize{$\cdots$}}
\put(87,81){\footnotesize{$\cdots$}}
\put(107,81){\footnotesize{$\cdots$}}
\put(124,81){\footnotesize{$\cdots$}}

\put(27,101){\footnotesize{$\xi_{0}$}}
\put(47,101){\footnotesize{$\xi_{1}$}}
\put(67,101){\footnotesize{$\xi_{2}$}}
\put(87,101){\footnotesize{$\cdots$}}
\put(107,101){\footnotesize{$\xi_{n}$}}
\put(124,101){\footnotesize{$\cdots$}}

\put(27,121){\footnotesize{$\xi_{0}$}}
\put(47,121){\footnotesize{$\xi_{1}$}}
\put(67,121){\footnotesize{$\xi_{2}$}}
\put(87,121){\footnotesize{$\cdots$}}
\put(107,121){\footnotesize{$\xi_{n}$}}
\put(124,121){\footnotesize{$\cdots$}}

\psline{->}(70,10)(90,10)
\put(79,6){$\rm{T}_1$}

\put(11,38){\footnotesize{$(0,1)$}}
\put(11,58){\footnotesize{$(0,2)$}}
\put(14,78){\footnotesize{$\vdots$}}
\put(11,98){\footnotesize{$(0,n)$}}
\put(4,118){\footnotesize{$(0,n+1)$}}

\psline{->}(10, 70)(10,90)
\put(5,80){$\rm{T}_2$}

\put(20,28){\footnotesize{$b$}}
\put(20,48){\footnotesize{$\eta_{1}$}}
\put(20,68){\footnotesize{$\eta_{2}$}}
\put(22,88){\footnotesize{$\vdots$}}
\put(20,108){\footnotesize{$\eta_{n}$}}
\put(22,128){\footnotesize{$\vdots$}}

\put(40,28){\footnotesize{$\eta_{0}$}}
\put(40,48){\footnotesize{$\eta_{1}$}}
\put(40,68){\footnotesize{$\eta_{2}$}}
\put(42,88){\footnotesize{$\vdots$}}
\put(40,108){\footnotesize{$\eta_{n}$}}
\put(42,128){\footnotesize{$\vdots$}}

\put(60,28){\footnotesize{$\eta_{0}$}}
\put(60,48){\footnotesize{$\eta_{1}$}}
\put(60,68){\footnotesize{$\eta_{2}$}}
\put(62,88){\footnotesize{$\vdots$}}
\put(60,108){\footnotesize{$\eta_{n}$}}
\put(62,128){\footnotesize{$\vdots$}}

\put(80,28){\footnotesize{$\eta_{0}$}}
\put(80,48){\footnotesize{$\eta_{1}$}}
\put(80,68){\footnotesize{$\eta_{2}$}}
\put(82,88){\footnotesize{$\vdots$}}
\put(80,108){\footnotesize{$\eta_{n}$}}
\put(82,128){\footnotesize{$\vdots$}}

\put(100,28){\footnotesize{$\eta_{0}$}}
\put(100,48){\footnotesize{$\eta_{1}$}}
\put(100,68){\footnotesize{$\eta_{2}$}}
\put(102,88){\footnotesize{$\vdots$}}
\put(100,108){\footnotesize{$\eta_{n}$}}
\put(102,128){\footnotesize{$\vdots$}}

\put(122,28){\footnotesize{$\vdots$}}
\put(122,48){\footnotesize{$\vdots$}}
\put(122,68){\footnotesize{$\vdots$}}
\put(122,88){\footnotesize{$\vdots$}}
\put(122,108){\footnotesize{$\vdots$}}
\put(122,128){\footnotesize{$\vdots$}}

\end{picture}
\end{center}
\caption{{}}
\label{Figure 1}
\end{figure}

\begin{proposition}
\label{propA}The $2$-variable weighted shift $\mathbf{T}$ defined by (\ref%
{eqnew1}) and (\ref{eqnew2}) is subnormal only if $a\leq s$, where $s:=\sqrt{%
\frac{\xi _{0}^{2}\xi _{1}^{2}\eta _{1}^{2}}{\xi _{1}^{2}\eta _{0}^{2}+\xi
_{0}^{2}\eta _{1}^{2}-\xi _{0}^{2}\eta _{0}^{2}}}$.

\begin{proof}
Suppose that $T$ above is subnormal, and let $\mu $ be the associated Berger
measure. \ Then the following partial moment matrix $M$, corresponding to
the moments of $\mu $ associated with the monomials $1,s,t$ and $ts$, must
be positive semi-definite:%
\begin{equation*}
M:=\left( 
\begin{array}{cccc}
1 & a^{2} & b^{2} & a^{2}\eta _{0}^{2} \\ 
a^{2} & a^{2}\xi _{1}^{2} & a^{2}\eta _{0}^{2} & a^{2}\eta _{0}^{2}\xi
_{1}^{2} \\ 
b^{2} & a^{2}\eta _{0}^{2} & b^{2}\eta _{1}^{2} & a^{2}\eta _{0}^{2}\eta
_{1}^{2} \\ 
a^{2}\eta _{0}^{2} & a^{2}\eta _{0}^{2}\xi _{1}^{2} & a^{2}\eta _{0}^{2}\eta
_{1}^{2} & a^{2}\eta _{0}^{2}\xi _{1}^{2}\eta _{1}^{2}%
\end{array}%
\right) .
\end{equation*}%
Now, using \textit{Mathematica} we obtain%
\begin{eqnarray*}
\det (M) &\geq &0 \\
&\Leftrightarrow &a^{6}\eta _{0}^{4}(\xi _{1}^{2}-\xi _{0}^{2})(\eta
_{1}^{2}-\eta _{0}^{2})(a^{2}\xi _{0}^{2}\eta _{0}^{2}-a^{2}\xi _{1}^{2}\eta
_{0}^{2}-a^{2}\xi _{0}^{2}\eta _{1}^{2}+\xi _{0}^{2}\xi _{1}^{2}\eta
_{1}^{2})\geq 0 \\
&\Leftrightarrow &a^{2}\xi _{0}^{2}\eta _{0}^{2}-a^{2}\xi _{1}^{2}\eta
_{0}^{2}-a^{2}\xi _{0}^{2}\eta _{1}^{2}+\xi _{0}^{2}\xi _{1}^{2}\eta
_{1}^{2}\geq 0 \\
&\Leftrightarrow &a\leq \sqrt{\frac{\xi _{0}^{2}\xi _{1}^{2}\eta _{1}^{2}}{%
\xi _{1}^{2}\eta _{0}^{2}+\xi _{0}^{2}\eta _{1}^{2}-\xi _{0}^{2}\eta _{0}^{2}%
}}=s.
\end{eqnarray*}
\end{proof}
\end{proposition}

\begin{proposition}
\label{propB}The $2$-variable weighted shift $\mathbf{T}$ defined by (\ref%
{eqnew1}) and (\ref{eqnew2}) is hyponormal if and only if $a\leq h$, where $%
h:=\xi _{0}\sqrt{\frac{\xi _{1}^{2}\eta _{1}^{2}-\xi _{0}^{2}\eta _{0}^{2}}{%
\xi _{0}^{2}\eta _{1}^{2}+\xi _{1}^{2}\eta _{0}^{2}-2\xi _{0}^{2}\eta
_{0}^{2}}}$.

\begin{proof}
From the definition of $\mathbf{T}$ and the Six-point Test (Theorem \ref%
{joint hypo}), it is clear that all we need is for the following matrix to
be positive semi-definite:%
\begin{equation*}
L:=\left( 
\begin{array}{cc}
\xi _{1}^{2}-a^{2} & \xi _{0}\eta _{0}-ab \\ 
\xi _{0}\eta _{0}-ab & \eta _{1}^{2}-b^{2}%
\end{array}%
\right) .
\end{equation*}%
Observe that 
\begin{eqnarray*}
\det L &\geq &0\Leftrightarrow \xi _{1}^{2}\eta _{1}^{2}-\xi
_{1}^{2}b^{2}-\xi _{0}^{2}\eta _{0}^{2}-a^{2}\eta _{1}^{2}+2ab\xi _{0}\eta
_{0}\geq 0 \\
&\Leftrightarrow &\xi _{1}^{2}\eta _{1}^{2}-\xi _{1}^{2}\frac{a^{2}\eta
_{0}^{2}}{\xi _{0}^{2}}-\xi _{0}^{2}\eta _{0}^{2}-a^{2}\eta
_{1}^{2}+2a^{2}\eta _{0}^{2}\geq 0\;\;\text{(using }b\xi _{0}=a\eta _{0}\;\;%
\text{(\ref{commuting2}))} \\
&\Leftrightarrow &a^{2}\leq \frac{\xi _{0}^{2}(\xi _{1}^{2}\eta _{1}^{2}-\xi
_{0}^{2}\eta _{0}^{2})}{\xi _{0}^{2}\eta _{1}^{2}+\xi _{1}^{2}\eta
_{0}^{2}-2\xi _{0}^{2}\eta _{0}^{2}}=h^{2}.
\end{eqnarray*}%
Thus, $a\leq h$ is clearly a necessary condition for the hyponormality of $%
\mathbf{T}$. \ Now, a straightforward calculation shows that $h<\xi _{1}$;
for, 
\begin{equation*}
\xi _{1}^{2}-h^{2}=\frac{\eta _{0}^{2}(\xi _{1}^{2}-\xi _{0}^{2})^{2}}{\xi
_{0}^{2}\eta _{1}^{2}+\xi _{1}^{2}\eta _{0}^{2}-2\xi _{0}^{2}\eta _{0}^{2}}%
>0.
\end{equation*}%
It follows that $a\leq h$ implies $a<\xi _{1}$, and therefore $L\geq 0$ by
the Nested Determinant Test \cite{Atk}. \ Thus, the condition $a\leq h$ is
also sufficient for the hyponormality of $\mathbf{T}$, and the proof is
complete.
\end{proof}
\end{proposition}

It follows from Propositions \ref{propA} and \ref{propB} that to ascertain
the existence of a nonsubnormal, hyponormal $2$-variable weighted shift $%
\mathbf{T}$ (with $T_{1}$ and $T_{2}$ subnormal), it suffices to show that
for appropriate choices of $\xi _{0},$ $\xi _{1},$ $\eta _{0}$ and $\eta
_{1} $, it is possible to obtain $s<h$, while keeping $a\leq \xi _{ext}(\nu
_{\mathcal{M}})$ and $b\equiv \frac{a\eta _{0}}{\xi _{0}}\leq \eta
_{ext}(\omega _{\mathcal{M}})$. \ Now, 
\begin{equation*}
h^{2}-s^{2}=\frac{\xi _{0}^{4}\eta _{0}^{2}(\xi _{1}^{2}-\xi _{0}^{2})(\eta
_{1}^{2}-\eta _{0}^{2})}{(\xi _{0}^{2}\eta _{1}^{2}+\xi _{1}^{2}\eta
_{0}^{2}-2\xi _{0}^{2}\eta _{0}^{2})(\xi _{1}^{2}\eta _{0}^{2}+\xi
_{0}^{2}\eta _{1}^{2}-\xi _{0}^{2}\eta _{0}^{2})}>0.
\end{equation*}%
Therefore, it suffices to prove the existence of strictly increasing weight
sequences $\{\xi _{i}\}$ and $\{\eta _{j}\}$ such that

\begin{enumerate}
\item $a\leq h\;$(hyponormality of $\mathbf{T}$)

\item $a>s\;$(nonsubnormality of $\mathbf{T}$)

\item $a\leq \xi _{ext}(\nu _{\mathcal{M}})$ (subnormality of $T_{1}$)

\item $a\leq s_{2}:=\frac{\xi _{0}}{\eta _{0}}\eta _{ext}(\omega _{\mathcal{M%
}})$ (subnormality of $T_{2}$).
\end{enumerate}

We now seek to determine the relative positions of $h$, $s$, $s_{2}$, $\xi
_{0}$, $\xi _{ext}(\nu _{\mathcal{M}})$ and $\xi _{1}$ in the positive real
axis.\newline
\textbf{Claim 1}: $\xi _{0}\leq \xi _{ext}(\nu _{\mathcal{M}})$. \ This is
straightforward from the fact that $shift\{\xi _{0},\xi _{1},...)$ is
subnormal. \ \newline
\textbf{Claim 2}: $\xi _{0}<s$. \ For, 
\begin{equation*}
s^{2}-\xi _{0}^{2}=\frac{\xi _{0}^{2}\xi _{1}^{2}\eta _{1}^{2}}{\xi
_{1}^{2}\eta _{0}^{2}+\xi _{0}^{2}\eta _{1}^{2}-\xi _{0}^{2}\eta _{0}^{2}}%
-\xi _{0}^{2}=\frac{\xi _{0}^{2}(\xi _{1}^{2}-\xi _{0}^{2})(\eta
_{1}^{2}-\eta _{0}^{2})}{\xi _{1}^{2}\eta _{0}^{2}+\xi _{0}^{2}\eta
_{1}^{2}-\xi _{0}^{2}\eta _{0}^{2}}>0.
\end{equation*}%
\newline
\textbf{Claim 3}: $s<\xi _{1}$. \ For,%
\begin{equation*}
\xi _{1}^{2}-s^{2}=\xi _{1}^{2}-\frac{\xi _{0}^{2}\xi _{1}^{2}\eta _{1}^{2}}{%
\xi _{1}^{2}\eta _{0}^{2}+\xi _{0}^{2}\eta _{1}^{2}-\xi _{0}^{2}\eta _{0}^{2}%
}=\frac{\xi _{1}^{2}\eta _{0}^{2}(\xi _{1}^{2}-\xi _{0}^{2})}{\xi
_{1}^{2}\eta _{0}^{2}+\xi _{0}^{2}\eta _{1}^{2}-\xi _{0}^{2}\eta _{0}^{2}}>0.
\end{equation*}%
\newline
\textbf{Claim 4}: $h<\xi _{1}$. \ This was established in the proof of
Proposition \ref{propB}.\newline
\textbf{Claim 5}: $s<s_{2}$ whenever $\eta _{0}<u:=\frac{\xi _{0}^{2}\eta
_{e}^{2}\eta _{1}^{2}}{\xi _{1}^{2}(\eta _{1}^{2}-\eta _{e}^{2})+\xi
_{0}^{2}\eta _{e}^{2}}$, where $\eta _{e}\equiv \eta _{ext}(\omega _{%
\mathcal{M}})$. \ For, 
\begin{eqnarray*}
s_{2}^{2}-s^{2} &=&\frac{\xi _{0}^{2}}{\eta _{0}^{2}}\eta _{e}^{2}-\frac{\xi
_{0}^{2}\xi _{1}^{2}\eta _{1}^{2}}{\xi _{1}^{2}\eta _{0}^{2}+\xi
_{0}^{2}\eta _{1}^{2}-\xi _{0}^{2}\eta _{0}^{2}} \\
&=&\xi _{0}^{2}\frac{\xi _{0}^{2}\eta _{e}^{2}\eta _{1}^{2}-\eta
_{0}^{2}[\xi _{1}^{2}(\eta _{1}^{2}-\eta _{e}^{2})+\xi _{0}^{2}\eta _{e}^{2}]%
}{\eta _{0}^{2}(\xi _{1}^{2}\eta _{0}^{2}+\xi _{0}^{2}\eta _{1}^{2}-\xi
_{0}^{2}\eta _{0}^{2})}.
\end{eqnarray*}%
It then follows that $s_{2}-s>0\Leftrightarrow \xi _{0}^{2}\eta _{e}^{2}\eta
_{1}^{2}-\eta _{0}^{2}[\xi _{1}^{2}(\eta _{1}^{2}-\eta _{e}^{2})+\xi
_{0}^{2}\eta _{e}^{2}]>0$, as desired.\newline
\textbf{Claim 6}: $h\leq s_{2}$ whenever $\eta _{0}\leq v:=\frac{\xi
_{1}^{2}(\eta _{1}^{2}-\eta _{e}^{2})+2\xi _{0}^{2}\eta _{e}^{2}-\sqrt{(\eta
_{1}^{2}-\eta _{e}^{2})(\xi _{1}^{4}(\eta _{1}^{2}-\eta _{e}^{2})+4\xi
_{0}^{2}\eta _{e}^{2}(\xi _{1}^{2}-\xi _{0}^{2}))}}{2\xi _{0}^{2}}$. \ Since 
\begin{equation*}
s_{2}^{2}-h^{2}=\frac{\xi _{0}^{2}\{\xi _{0}^{2}\eta _{0}^{4}-[\xi
_{1}^{2}(\eta _{1}^{2}-\eta _{e}^{2})+2\xi _{0}^{2}\eta _{e}^{2}]\eta
_{0}^{2}+\xi _{0}^{2}\eta _{e}^{2}\eta _{1}^{2}\}}{\eta _{0}^{2}(\xi
_{0}^{2}\eta _{1}^{2}+\xi _{1}^{2}\eta _{0}^{2}-2\xi _{0}^{2}\eta _{0}^{2})},
\end{equation*}%
it follows that $h\leq s_{2}$ if and only if the quadratic form 
\begin{eqnarray*}
q(t) &\equiv &At^{2}+Bt+C \\
&:&=\xi _{0}^{2}t^{2}-[\xi _{1}^{2}(\eta _{1}^{2}-\eta _{e}^{2})+2\xi
_{0}^{2}\eta _{e}^{2}]t+\xi _{0}^{2}\eta _{e}^{2}\eta _{1}^{2}
\end{eqnarray*}%
is nonnegative. \ Since $A$ and $C$ are positive, and $B$ is negative, we
need to study the discriminant, $\Delta :=B^{2}-4AC$. \ Now, 
\begin{eqnarray*}
\Delta &=&(\xi _{1}^{2}(\eta _{1}^{2}-\eta _{e}^{2})+2\xi _{0}^{2}\eta
_{e}^{2})^{2}-4\xi _{0}^{4}\eta _{e}^{2}\eta _{1}^{2} \\
&=&(\eta _{1}^{2}-\eta _{e}^{2})[\xi _{1}^{4}\eta _{1}^{2}-\eta
_{e}^{2}(2\xi _{0}^{2}-\xi _{1}^{2})^{2}],
\end{eqnarray*}%
so $\Delta \geq 0\Leftrightarrow \xi _{1}^{4}\eta _{1}^{2}-\eta
_{e}^{2}(2\xi _{0}^{2}-\xi _{1}^{2})^{2}\geq 0$. \ Since $\xi _{1}^{4}\eta
_{1}^{2}-\eta _{e}^{2}(2\xi _{0}^{2}-\xi _{1}^{2})^{2}=\xi _{1}^{4}(\eta
_{1}^{2}-\eta _{e}^{2})+4\xi _{0}^{2}\eta _{e}^{2}(\xi _{1}^{2}-\xi
_{0}^{2}) $, we see that $\Delta $ is always positive. \ We conclude that $%
q\geq 0$ on the interval $[0,t_{1}]$, where $t_{1}:=\frac{-B-\sqrt{\Delta }}{%
2A}$ is the leftmost zero of $q$. \ Finally, a straightforward calculation
shows that $t_{1}=v$.

We now summarize what we have so far. \ For $\eta _{0}<\min \{u,v\}$ we have 
\begin{equation*}
\left\{ 
\begin{array}{c}
\xi _{0}<s<h\leq s_{2} \\ 
\\ 
h<\xi _{1} \\ 
\\ 
\xi _{ext}(\nu _{\mathcal{M}})<\xi _{1}.%
\end{array}%
\right.
\end{equation*}%
Thus, if we can ensure that $h\leq \xi _{ext}(\nu _{\mathcal{M}})$, the
construction of the example will be complete by taking $a$ such that $%
s<a\leq h$. \ Now, since $h\leq s_{2}$, an easy way to accomplish this is to
build $shift(\xi _{0},\xi _{1},...)$ in such a way that $\xi _{ext}(\nu _{%
\mathcal{M}})=s_{2}$. \ To do this, we appeal to Lemma \ref{abclem}, that
is, we first build a $2$-step recursively generated weighted shift whose
first three weights are $s_{2}$, $\xi _{1}$ and $\xi _{2}$, and we then
consider the shift $W_{\xi _{0}(\xi _{1},\xi _{2},\xi _{3})\symbol{94}}$,
where $\xi _{3}$ is given by $\xi _{3}:=\frac{\varphi _{0}}{\xi _{2}^{2}}%
+\varphi _{1}$ obtained from the equation $\gamma _{4}=\varphi _{0}\gamma
_{2}+\varphi _{1}\gamma _{3}$. \ Observe that the extremal value of $W_{(\xi
_{1},\xi _{2},\xi _{3})\symbol{94}}$ is $s_{2}$, and that $\xi _{0}<s_{2}$,
so the subnormality of $W_{\xi _{0}(\xi _{1},\xi _{2},\xi _{3})\symbol{94}}$
is guaranteed. \ This completes the construction of the example.

\begin{theorem}
Let $\mathbf{T\equiv (}T_{1},T_{2})$ be the $2$-variable weighted shift
defined by (\ref{eqnew1}) and (\ref{eqnew2}), let 
\begin{equation*}
\left\{ 
\begin{array}{l}
h:=\xi _{0}\sqrt{\frac{\xi _{1}^{2}\eta _{1}^{2}-\xi _{0}^{2}\eta _{0}^{2}}{%
\xi _{0}^{2}\eta _{1}^{2}+\xi _{1}^{2}\eta _{0}^{2}-2\xi _{0}^{2}\eta
_{0}^{2}}}, \\ 
\\ 
s:=\sqrt{\frac{\xi _{0}^{2}\xi _{1}^{2}\eta _{1}^{2}}{\xi _{1}^{2}\eta
_{0}^{2}+\xi _{0}^{2}\eta _{1}^{2}-\xi _{0}^{2}\eta _{0}^{2}}}, \\ 
\\ 
s_{2}:=\frac{\xi _{0}}{\eta _{0}}\eta _{e}\text{, where }\eta _{e}\equiv
\eta _{ext}(\omega _{\mathcal{M}}), \\ 
\\ 
u:=\frac{\xi _{0}^{2}\eta _{e}^{2}\eta _{1}^{2}}{\xi _{1}^{2}(\eta
_{1}^{2}-\eta _{e}^{2})+\xi _{0}^{2}\eta _{e}^{2}}\text{, and } \\ 
\\ 
v:=\frac{\xi _{1}^{2}(\eta _{1}^{2}-\eta _{e}^{2})+2\xi _{0}^{2}\eta
_{e}^{2}-\sqrt{(\eta _{1}^{2}-\eta _{e}^{2})(\xi _{1}^{4}(\eta _{1}^{2}-\eta
_{e}^{2})+4\xi _{0}^{2}\eta _{e}^{2}(\xi _{1}^{2}-\xi _{0}^{2}))}}{2\xi
_{0}^{2}}.%
\end{array}%
\right.
\end{equation*}%
Assume further that, as above, $s_{2}=\xi _{ext}(\nu _{\mathcal{M}})$ and $%
\eta _{0}\leq \min \{u,v\}$. \ Finally, choose $a$ such that $s<a\leq h$. \
Then

\begin{enumerate}
\item[(i)] $T_{1}T_{2}=T_{2}T_{1}$;\medskip

\item[(ii)] $T_{1}$ is subnormal;\medskip

\item[(iii)] $T_{2}$ is subnormal;\medskip

\item[(iv)] $\mathbf{T}$ is hyponormal; and\medskip

\item[(v)] $\mathbf{T}$ is not subnormal.
\end{enumerate}
\end{theorem}

\begin{example}
\label{counterex1}For a concrete numerical example, let $d\omega _{\mathcal{M%
}}(t):=2dt$ on $[\frac{1}{2},1]$, so that $\left\| \frac{1}{t}\right\|
_{L^{1}(\omega _{\mathcal{M}})}=2\ln 2$. \ It follows that $\eta _{e}\equiv
\eta _{ext}(\omega _{\mathcal{M}})=\frac{1}{\sqrt{2\ln 2}}$ and $\eta _{1}=%
\frac{\sqrt{3}}{2}$. \ Now take $\xi _{0}:=\frac{1}{2}$ and $\xi _{1}:=1$. \
Then $u=\frac{1}{4\left( 2\ln 2-1\right) }\cong 0.647$ and $v=\frac{1}{4}%
\frac{6\ln 2-2-\sqrt{2}\sqrt{\left( 3\ln 2-2\right) }\sqrt{\left( 6\ln
2-1\right) }}{\ln 2}\cong 0.523$, so we can take $\eta _{0}:=\frac{1}{2}$. \
With this choice of $\eta _{0}$ we obtain $s=\frac{\sqrt{2}}{2}\cong 0.707$, 
$h=\frac{1}{2}\sqrt{\frac{11}{5}}\cong 0.742$ and $s_{2}=\eta _{e}=\frac{1}{%
\sqrt{2\ln 2}}\cong 0.849$. \ We can then take $a\in (s,h]$, for instance $%
a:=0.72$. \ To build the weighted shift $W_{\xi }$ we start with $s_{2}$, $%
\xi _{1}$ and $\xi _{2}:=\sqrt{2}$ to obtain $\varphi _{0}=\frac{1}{1-2\ln 2}
$ and $\varphi _{1}=\frac{1-4\ln 2}{1-2\ln 2}$. \ This gives $\xi _{3}=\frac{%
1}{2}\sqrt{\frac{16\ln 2-5}{2\ln 2-1}}\cong 1.985$. \ The $2$-atomic measure 
$\nu _{\mathcal{M}}$ for $W_{(\xi _{1},\xi _{2},\xi _{3})\symbol{94}}$ has
atoms $t_{0}\cong 0.659$ and $t_{1}\cong 3.93$ and densities $\rho _{0}\cong
0.981$ and $\rho _{1}\cong 0.019$. \ With these values we can compute $%
\left\| \frac{1}{t}\right\| _{L^{1}(\nu _{\mathcal{M}})}\cong 1.494$;
observe that $\xi _{0}=\frac{1}{2}\leq \xi _{ext}(\nu _{\mathcal{M}})=\sqrt{%
\frac{1}{\left\| \frac{1}{t}\right\| _{L^{1}(\nu _{\mathcal{M}})}}}\cong
0.818$. \ By Proposition \ref{backward}, the measure associated to $%
shift(\xi _{0},\xi _{1},\xi _{2},...)$ is $d\nu (t)=\frac{1}{4t}(\rho
_{0}d\delta _{t_{0}}(t)+\rho _{1}d\delta _{t_{1}}(t))+(1-\frac{1}{4}\left\| 
\frac{1}{t}\right\| _{L^{1}(\nu _{\mathcal{M}})})d\delta _{0}(t)$.
\end{example}

\section{The Second Family of Counterexamples}

Recall that a unilateral weighted shift $W_{\alpha }$ is subnormal if and
only if there exists a probability measure $\xi \equiv \xi _{\alpha }$
supported in $[0,\left\| W_{\alpha }\right\| ^{2}]$ such that $\gamma
_{k}(\alpha ):=\alpha _{0}^{2}\cdot ...\cdot \alpha _{k-1}^{2}=\int
t^{k}\;d\xi (t)\;\;(k\geq 1)$. \ For instance, when $\alpha _{1}=\alpha
_{2}=...=1$ (i.e., $W_{\alpha }\equiv shift(\alpha _{0},1,1,...)$), we have $%
\xi _{\alpha }=(1-\alpha _{0}^{2})\delta _{0}+\alpha _{0}^{2}\delta _{1}$. \
The proof of the following lemma is straightforward.

\begin{lemma}
\label{lemtensor}Given two $1$-variable weight sequences $\alpha $ and $%
\beta $, the $2$-variable weighted shift \newline
$(W_{\alpha }\bigotimes I,I\bigotimes W_{\beta })$ is always subnormal, with
Berger measure $\mu :=\xi _{\alpha }\times \xi _{\beta }$.
\end{lemma}

\begin{definition}
Let $\mu $ and $\nu $ be two positive measures on $\mathbb{R}_{+}.$ \ We say
that $\mu \leq \nu $ on $X:=\mathbb{R}_{+},$ if $\mu (E)\leq \nu (E)$ for
all Borel subset $E\subseteq \mathbb{R}_{+}$; equivalently, $\mu \leq \nu $
if and only if $\int fd\mu \leq \int fd\nu $ for all $f\in C(X)$ such that $%
f\geq 0$ on $\mathbb{R}_{+}$.
\end{definition}

\begin{definition}
Let $\mu $ be a probability measure on $X\times Y$, and assume that $\frac{1%
}{t}\in L^{1}(\mu ).$ \ The \emph{extremal measure} $\mu _{ext}$ (which is
also a probability measure) on $X\times Y$ is given by $d\mu
_{ext}(s,t):=(1-\delta _{0}(t))\frac{1}{t\left\| \frac{1}{t}\right\|
_{L^{1}(\mu )}}d\mu (s,t)$. \ 
\end{definition}

\begin{example}
Let $B_{+}$ be the Bergman shift on $\ell ^{2}(\mathbb{Z}_{+})$ and let $%
\mathcal{M}\equiv \mathcal{M}_{1}:=\bigvee \{e_{1},e_{2},...\}$. \ The shift 
$B_{+}|_{\mathcal{M}}$ is subnormal, with Berger measure $d\mu (t):=tdt$ on $%
[0,1]$. \ Then $d\mu _{ext}(t)=dt$, so the extremal measure $\mu _{ext}$ is
the Berger measure of $B_{+}$.
\end{example}

\begin{definition}
\label{defmarg}Given a measure $\mu $ on $X\times Y$, the marginal measure $%
\mu ^{X}$ is given by $\mu ^{X}:=\mu \circ \pi _{X}^{-1}$, where $\pi
_{X}:X\times Y\rightarrow X$ is the canonical projection onto $X$. \ Thus, $%
\mu ^{X}(E)=\mu (E\times Y)$, for every $E\subseteq X$. \ Observe that if $%
\mu $ is a probability measure, then so is $\mu ^{X}$.
\end{definition}

\begin{lemma}
\label{lemBerger}Let $\mu $ be the Berger measure of a $2$-variable weighted
shift $\mathbf{T}$ and let $\nu $ be the Berger measure of $shift(\alpha
_{00},\alpha _{10},...)$. \ Then $\nu =\mu ^{X}$. \ As a consequence, $\iint
f(s)\;d\mu (s,t)=\int f(s)\;d\mu ^{X}(s)$ for all $f\in C(X)$.
\end{lemma}

\begin{proof}
Observe that $\int s^{i}\;d\nu (s)=\gamma _{i0}=\iint s^{i}\;d\mu (s,t)$ for
all $i\geq 0$. \ It follows that $\int f(s)\;d\nu (s)=\iint f(s)\;d\mu (s,t)$
for all $f\in C(X)$. \ Then, for any Borel set $E\subseteq X$, we have 
\begin{equation*}
\nu (E)=\int \chi _{E}\;d\nu =\iint \chi _{E\times Y}\;d\mu =\mu (E\times
Y)=\mu ^{X}(E),
\end{equation*}%
as desired. \ The second assertion follows immediately from what we have
established.
\end{proof}

\begin{corollary}
Let $\mu $ be the Berger measure of a $2$-variable weighted shift $\mathbf{T}
$. \ For $j\geq 1$, let $d\mu _{j}(s,t):=\frac{1}{\gamma _{0j}}t^{j}d\mu
(s,t)$. \ Then the Berger measure of $shift(\alpha _{0j},\alpha _{1j},...)$
is $\nu _{j}\equiv \mu _{j}^{X}$.\ 
\end{corollary}

\begin{example}
Let $\mu :=\xi \times \eta $ be a probability product measure on $X\times Y$%
. \ Then $\mu ^{X}=\xi $.
\end{example}

\begin{lemma}
\label{lemmarg}Let $\mu $ and $\omega $ be two measures on $X\times Y$, and
assume that $\mu \leq \omega $. \ Then $\mu ^{X}\leq \omega ^{X}$.
\end{lemma}

\begin{proof}
Straightforward from Definition \ref{defmarg}.
\end{proof}

\begin{proposition}
\label{backext}(Subnormal backward extension of a $2$-variable weighted
shift) \ Consider the following $2$-variable weighted shift (see Figure \ref%
{Figure 2}), and let $\mathcal{M}$ be the subspace associated to indices $%
\mathbf{k}$ with $k_{2}\geq 1$. \ Assume that $\mathbf{T}_{\mathcal{M}}$ is
subnormal with associated measure $\mu _{\mathcal{M}}$ and that $%
W_{0}:=shift(\alpha _{00},\alpha _{10},\cdots )$ is subnormal with
associated measure $\nu $. \ Then $\mathbf{T}$ is subnormal if and only if%
\newline
(i) $\ \frac{1}{t}\in L^{1}(\mu _{\mathcal{M}})$;\newline
(ii) $\ \beta _{00}^{2}\leq (\left\| \frac{1}{t}\right\| _{L^{1}(\mu _{%
\mathcal{M}})})^{-1}$;\newline
(iii) $\ \beta _{00}^{2}\left\| \frac{1}{t}\right\| _{L^{1}(\mu _{\mathcal{M}%
})}(\mu _{\mathcal{M}})_{ext}^{X}\leq \nu $.\newline
Moreover, if $\beta _{00}^{2}\left\| \frac{1}{t}\right\| _{L^{1}(\mu _{%
\mathcal{M}})}=1,$ then $(\mu _{\mathcal{M}})_{ext}^{X}=\nu $. \ In the case
when $\mathbf{T}$ is subnormal, the Berger measure $\mu $ of $\mathbf{T}$ is
given by 
\begin{equation*}
d\mu (s,t)=\beta _{00}^{2}\left\| \frac{1}{t}\right\| _{L^{1}(\mu _{\mathcal{%
M}})}d(\mu _{\mathcal{M}})_{ext}(s,t)+(d\nu (s)-\beta _{00}^{2}\left\| \frac{%
1}{t}\right\| _{L^{1}(\mu _{\mathcal{M}})}d(\mu _{\mathcal{M}%
})_{ext}^{X}(s))d\delta _{0}(t).
\end{equation*}%
\setlength{\unitlength}{1mm} \psset{unit=1mm} 
\begin{figure}[th]
\begin{center}
\begin{picture}(140,138)

\psline{->}(20,20)(135,20)
\psline(20,40)(125,40)
\psline(20,60)(125,60)
\psline(20,80)(125,80)
\psline(20,100)(125,100)
\psline(20,120)(125,120)
\psline{->}(20,20)(20,135)
\psline(40,20)(40,125)
\psline(60,20)(60,125)
\psline(80,20)(80,125)
\psline(100,20)(100,125)
\psline(120,20)(120,125)

\put(11,16){\footnotesize{$(0,0)$}}
\put(35,16){\footnotesize{$(1,0)$}}
\put(55,16){\footnotesize{$(2,0)$}}
\put(78,16){\footnotesize{$\cdots$}}
\put(95,16){\footnotesize{$(n,0)$}}
\put(115,16){\footnotesize{$(n+1,0)$}}

\put(27,21){\footnotesize{$\alpha_{0,0}$}}
\put(47,21){\footnotesize{$\alpha_{1,0}$}}
\put(67,21){\footnotesize{$\alpha_{2,0}$}}
\put(87,21){\footnotesize{$\cdots$}}
\put(107,21){\footnotesize{$\alpha_{n,0}$}}
\put(124,21){\footnotesize{$\alpha_{n+1,0}$}}

\put(27,41){\footnotesize{$\alpha_{0,1}$}}
\put(47,41){\footnotesize{$\alpha_{1,1}$}}
\put(67,41){\footnotesize{$\alpha_{2,1}$}}
\put(87,41){\footnotesize{$\cdots$}}
\put(107,41){\footnotesize{$\alpha_{n,1}$}}
\put(124,41){\footnotesize{$\cdots$}}

\put(27,61){\footnotesize{$\alpha_{0,2}$}}
\put(47,61){\footnotesize{$\alpha_{1,2}$}}
\put(67,61){\footnotesize{$\alpha_{2,2}$}}
\put(87,61){\footnotesize{$\cdots$}}
\put(107,61){\footnotesize{$\alpha_{n,2}$}}
\put(124,61){\footnotesize{$\cdots$}}

\put(27,81){\footnotesize{$\cdots$}}
\put(47,81){\footnotesize{$\cdots$}}
\put(67,81){\footnotesize{$\cdots$}}
\put(87,81){\footnotesize{$\cdots$}}
\put(107,81){\footnotesize{$\cdots$}}
\put(124,81){\footnotesize{$\cdots$}}

\put(27,101){\footnotesize{$\alpha_{0,n}$}}
\put(47,101){\footnotesize{$\alpha_{1,n}$}}
\put(67,101){\footnotesize{$\alpha_{2,n}$}}
\put(87,101){\footnotesize{$\cdots$}}
\put(107,101){\footnotesize{$\alpha_{n,n}$}}
\put(124,101){\footnotesize{$\cdots$}}

\put(27,121){\footnotesize{$\alpha_{0,n+1}$}}
\put(47,121){\footnotesize{$\alpha_{1,n+1}$}}
\put(67,121){\footnotesize{$\alpha_{2,n+1}$}}
\put(87,121){\footnotesize{$\cdots$}}
\put(107,121){\footnotesize{$\alpha_{n,n+1}$}}
\put(124,121){\footnotesize{$\cdots$}}

\psline{->}(70,14)(90,14)
\put(79,10){$\rm{T}_1$}

\put(11,38){\footnotesize{$(0,1)$}}
\put(11,58){\footnotesize{$(0,2)$}}
\put(14,78){\footnotesize{$\vdots$}}
\put(11,98){\footnotesize{$(0,n)$}}
\put(4,118){\footnotesize{$(0,n+1)$}}

\psline{->}(13,70)(13,90)
\put(5,80){$\rm{T}_2$}

\put(20,28){\footnotesize{$\beta_{0,0}$}}
\put(20,48){\footnotesize{$\beta_{0,1}$}}
\put(20,68){\footnotesize{$\beta_{0,2}$}}
\put(22,88){\footnotesize{$\vdots$}}
\put(20,108){\footnotesize{$\beta_{0,n}$}}
\put(22,128){\footnotesize{$\beta_{0,n+1}$}}

\put(40,28){\footnotesize{$\sqrt{\frac{\gamma_{1,1}}{\gamma_{1,0}}}$}}
\put(40,48){\footnotesize{$\sqrt{\frac{\gamma_{1,2}}{\gamma_{1,1}}}$}}
\put(40,68){\footnotesize{$\sqrt{\frac{\gamma_{1,3}}{\gamma_{1,2}}}$}}
\put(42,88){\footnotesize{$\vdots$}}
\put(40,108){\footnotesize{$\sqrt{\frac{\gamma_{1,n+1}}{\gamma_{1,n}}}$}}
\put(42,128){\footnotesize{$\vdots$}}

\put(60,28){\footnotesize{$\sqrt{\frac{\gamma_{2,1}}{\gamma_{2,0}}}$}}
\put(60,48){\footnotesize{$\sqrt{\frac{\gamma_{2,2}}{\gamma_{2,1}}}$}}
\put(60,68){\footnotesize{$\sqrt{\frac{\gamma_{2,3}}{\gamma_{2,2}}}$}}
\put(62,88){\footnotesize{$\vdots$}}
\put(60,108){\footnotesize{$\sqrt{\frac{\gamma_{2,n+1}}{\gamma_{2,n}}}$}}
\put(62,128){\footnotesize{$\vdots$}}

\put(82,28){\footnotesize{$\vdots$}}
\put(82,48){\footnotesize{$\vdots$}}
\put(82,68){\footnotesize{$\vdots$}}
\put(82,88){\footnotesize{$\vdots$}}
\put(82,128){\footnotesize{$\vdots$}}

\put(100,28){\footnotesize{$\sqrt{\frac{\gamma_{n,1}}{\gamma_{n,0}}}$}}
\put(100,48){\footnotesize{$\sqrt{\frac{\gamma_{n,2}}{\gamma_{n,1}}}$}}
\put(100,68){\footnotesize{$\sqrt{\frac{\gamma_{n,3}}{\gamma_{n,2}}}$}}
\put(102,88){\footnotesize{$\vdots$}}
\put(100,108){\footnotesize{$\sqrt{\frac{\gamma_{n,n+1}}{\gamma_{n,n}}}$}}
\put(102,128){\footnotesize{$\vdots$}}

\put(122,28){\footnotesize{$\vdots$}}
\put(122,48){\footnotesize{$\vdots$}}
\put(122,68){\footnotesize{$\vdots$}}
\put(122,88){\footnotesize{$\vdots$}}
\put(122,128){\footnotesize{$\vdots$}}

\end{picture}
\end{center}
\caption{Weight diagram of the 2-variable weighted shift in Proposition \ref%
{backext}}
\label{Figure 2}
\end{figure}
\end{proposition}

\begin{proof}
($\Rightarrow $) \ First, observe that the moments of $\mathbf{T}$ and $%
\mathbf{T}_{\mathcal{M}}$ are related as follows:%
\begin{equation}
\gamma _{\mathbf{k}+\mathbf{\varepsilon }_{2}}(\mathbf{T})=\beta
_{00}^{2}\gamma _{\mathbf{k}}(\mathbf{T}_{\mathcal{M}})\;(\text{all }\mathbf{%
k}\in \mathbb{Z}_{+}^{2}),  \label{momrel}
\end{equation}%
so under the assumption that $\mathbf{T}$ is subnormal we must have 
\begin{equation*}
\iint s^{i}t^{j}(td\mu )(s,t)=\iint s^{i}t^{j+1}d\mu (s,t)=\gamma _{i,j+1}(%
\mathbf{T})=\beta _{00}^{2}\gamma _{ij}=\iint s^{i}t^{j}\;d\mu _{\mathcal{M}%
}(s,t).
\end{equation*}%
Thus $td\mu (s,t)=\beta _{00}^{2}d\mu _{\mathcal{M}}(s,t)$ and $\mu _{%
\mathcal{M}}(E\times \{0\})=0$ for all $E\subseteq X$. \ It follows at once
that 
\begin{eqnarray*}
\iint \frac{1}{t}\;d\mu _{\mathcal{M}}(s,t) &=&\iint_{(t>0)}\frac{1}{t}%
\;d\mu _{\mathcal{M}}(s,t)=\frac{1}{\beta _{00}^{2}}\iint_{(t>0)}\frac{1}{t}%
\;td\mu (s,t) \\
&=&\frac{1}{\beta _{00}^{2}}\mu ((t>0))\leq \frac{1}{\beta _{00}^{2}},
\end{eqnarray*}%
which establishes parts (i) and (ii). \ As for part (iii), let $E\subseteq X$
and $F\subseteq Y$ be two arbitrary Borel sets. \ Then 
\begin{eqnarray}
\beta _{00}^{2}\left\| \frac{1}{t}\right\| _{L^{1}(\mu _{\mathcal{M}})}(\mu
_{\mathcal{M}})_{ext}(E\times F) &=&\beta _{00}^{2}\left\| \frac{1}{t}%
\right\| _{L^{1}(\mu _{\mathcal{M}})}\iint_{E\times F}(1-\delta _{0}(t))%
\frac{1}{t\left\| \frac{1}{t}\right\| _{L^{1}(\mu )}}d\mu _{\mathcal{M}}(s,t)
\notag \\
&=&\iint_{E\times (F\;\backslash \;\{0\})}\frac{1}{t}\beta _{00}^{2}d\mu _{%
\mathcal{M}}(s,t)=\mu _{\mathcal{M}}(E\times (F\;\backslash \;\{0\}))  \notag
\\
&\leq &\mu (E\times F),  \label{b00}
\end{eqnarray}%
and by Lemmas \ref{lemmarg} and \ref{lemBerger}, $\beta _{00}^{2}\left\| 
\frac{1}{t}\right\| _{L^{1}(\mu _{\mathcal{M}})}(\mu _{\mathcal{M}%
})_{ext}^{X}\leq \mu ^{X}=\nu $. \ Finally, observe that when $\beta
_{00}^{2}\left\| \frac{1}{t}\right\| _{L^{1}(\mu _{\mathcal{M}})}=1$, the
inequality in (\ref{b00}) becomes an equality, and therefore $(\mu _{%
\mathcal{M}})_{ext}^{X}=\nu $.

($\Leftarrow $) Assume that (i), (ii) and (iii) hold, and let 
\begin{equation*}
\mu :=\beta _{00}^{2}\left\| \frac{1}{t}\right\| _{L^{1}(\mu _{\mathcal{M}%
})}(\mu _{\mathcal{M}})_{ext}+[\nu -\beta _{00}^{2}\left\| \frac{1}{t}%
\right\| _{L^{1}(\mu _{\mathcal{M}})}(\mu _{\mathcal{M}})_{ext}^{X}]\times
\delta _{0}.\ 
\end{equation*}%
Of course, if $\beta _{00}^{2}\left\| \frac{1}{t}\right\| _{L^{1}(\mu _{%
\mathcal{M}})}=1$, then $\mu :=(\mu _{\mathcal{M}})_{ext}$, since the total
mass of the second summand is zero. \ We now compute the moments of $\mu $
and verify that they agree with the moments of $\mathbf{T}$. \ If $j>0$,
then 
\begin{eqnarray*}
\iint s^{i}t^{j}\;d\mu (s,t) &=&\beta _{00}^{2}\left\| \frac{1}{t}\right\|
_{L^{1}(\mu _{\mathcal{M}})}\iint s^{i}t^{j}\;d(\mu _{\mathcal{M}%
})_{ext}(s,t) \\
&=&\beta _{00}^{2}\left\| \frac{1}{t}\right\| _{L^{1}(\mu _{\mathcal{M}%
})}\iint s^{i}t^{j}(1-\delta _{0}(t))\frac{1}{t\left\| \frac{1}{t}\right\|
_{L^{1}(\mu )}}d\mu _{\mathcal{M}}(s,t) \\
&=&\beta _{00}^{2}\iint s^{i}t^{j-1}\;d\mu _{\mathcal{M}}(s,t)=\beta
_{00}^{2}\gamma _{(i,j-1)}(\mathbf{T}_{\mathcal{M}})=\gamma _{(i,j)}(\mathbf{%
T})\text{,}
\end{eqnarray*}%
as desired. \ When $j=0$, we have 
\begin{eqnarray*}
\iint s^{i}\;d\mu (s,t) &=&\beta _{00}^{2}\left\| \frac{1}{t}\right\|
_{L^{1}(\mu _{\mathcal{M}})}\iint s^{i}\;d(\mu _{\mathcal{M}})_{ext}(s,t) \\
&&+\int s^{i}\;d(\nu -\beta _{00}^{2}\left\| \frac{1}{t}\right\| _{L^{1}(\mu
_{\mathcal{M}})}(\mu _{\mathcal{M}})_{ext}^{X})(s) \\
&=&\beta _{00}^{2}\left\| \frac{1}{t}\right\| _{L^{1}(\mu _{\mathcal{M}%
})}\int s^{i}\;d(\mu _{\mathcal{M}})_{ext}^{X}(s) \\
&&+\int s^{i}\;d\nu (s)-\beta _{00}^{2}\left\| \frac{1}{t}\right\|
_{L^{1}(\mu _{\mathcal{M}})}\int s^{i}\;d(\mu _{\mathcal{M}})_{ext}^{X}(s) \\
&&\text{(using Lemma \ref{lemBerger} for the first term)} \\
&=&\int s^{i}\;d\nu (s)=\gamma _{(i,0)}(\mathbf{T})\text{,}
\end{eqnarray*}%
as desired. \ It follows that $\mathbf{T}$ is subnormal, with Berger measure 
$\mu $.
\end{proof}

We are now ready to exhibit our second family of counterexamples to
Conjecture \ref{conjecture}. \ Consider the $2$-variable weighted shift
given by Figure \ref{Figure 3}, where $\max \{y,x,\frac{ay}{x}\}<1$.

\setlength{\unitlength}{1mm} \psset{unit=1mm} 
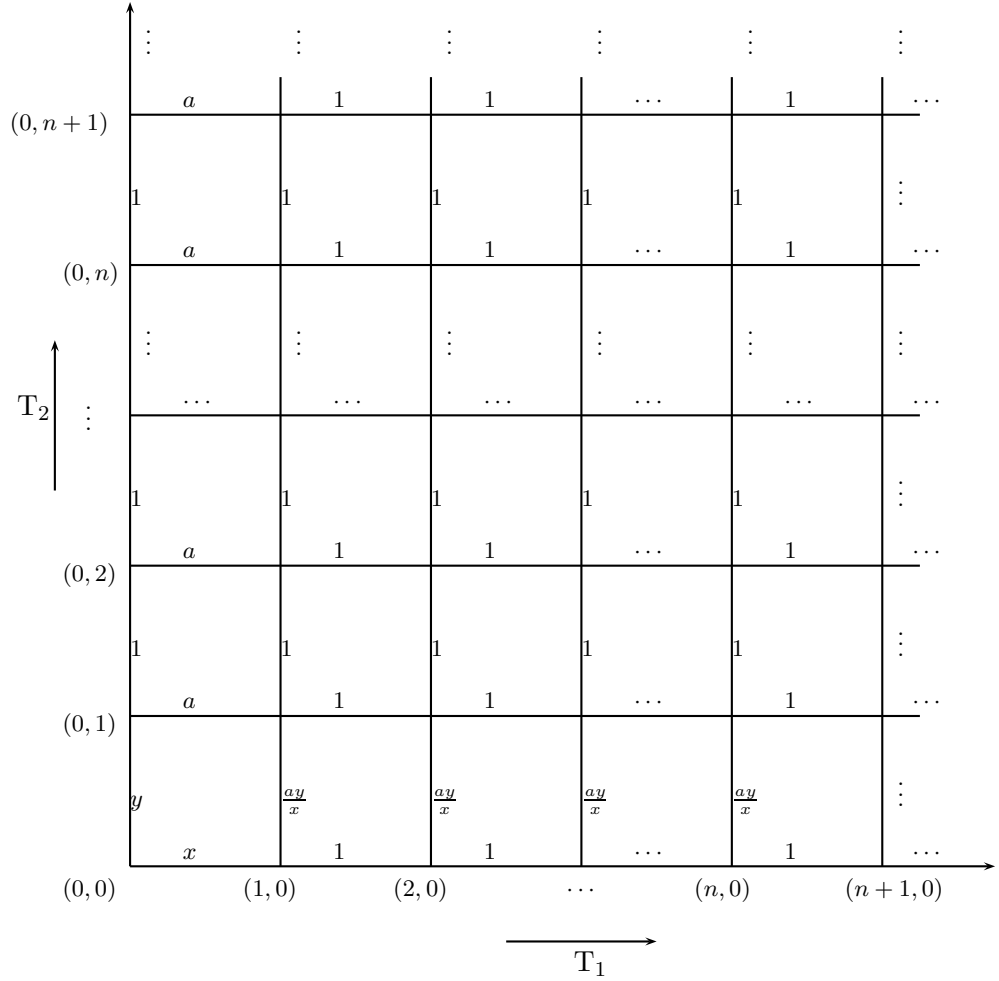
\begin{figure}[th]
\begin{center}
\begin{picture}(140,138)

\psline{->}(20,20)(135,20)
\psline(20,40)(125,40)
\psline(20,60)(125,60)
\psline(20,80)(125,80)
\psline(20,100)(125,100)
\psline(20,120)(125,120)
\psline{->}(20,20)(20,135)
\psline(40,20)(40,125)
\psline(60,20)(60,125)
\psline(80,20)(80,125)
\psline(100,20)(100,125)
\psline(120,20)(120,125)

\put(11,16){\footnotesize{$(0,0)$}}
\put(35,16){\footnotesize{$(1,0)$}}
\put(55,16){\footnotesize{$(2,0)$}}
\put(78,16){\footnotesize{$\cdots$}}
\put(95,16){\footnotesize{$(n,0)$}}
\put(115,16){\footnotesize{$(n+1,0)$}}

\put(27,21){\footnotesize{$x$}}
\put(47,21){\footnotesize{$1$}}
\put(67,21){\footnotesize{$1$}}
\put(87,21){\footnotesize{$\cdots$}}
\put(107,21){\footnotesize{$1$}}
\put(124,21){\footnotesize{$\cdots$}}

\put(27,41){\footnotesize{$a$}}
\put(47,41){\footnotesize{$1$}}
\put(67,41){\footnotesize{$1$}}
\put(87,41){\footnotesize{$\cdots$}}
\put(107,41){\footnotesize{$1$}}
\put(124,41){\footnotesize{$\cdots$}}

\put(27,61){\footnotesize{$a$}}
\put(47,61){\footnotesize{$1$}}
\put(67,61){\footnotesize{$1$}}
\put(87,61){\footnotesize{$\cdots$}}
\put(107,61){\footnotesize{$1$}}
\put(124,61){\footnotesize{$\cdots$}}

\put(27,81){\footnotesize{$\cdots$}}
\put(47,81){\footnotesize{$\cdots$}}
\put(67,81){\footnotesize{$\cdots$}}
\put(87,81){\footnotesize{$\cdots$}}
\put(107,81){\footnotesize{$\cdots$}}
\put(124,81){\footnotesize{$\cdots$}}

\put(27,101){\footnotesize{$a$}}
\put(47,101){\footnotesize{$1$}}
\put(67,101){\footnotesize{$1$}}
\put(87,101){\footnotesize{$\cdots$}}
\put(107,101){\footnotesize{$1$}}
\put(124,101){\footnotesize{$\cdots$}}

\put(27,121){\footnotesize{$a$}}
\put(47,121){\footnotesize{$1$}}
\put(67,121){\footnotesize{$1$}}
\put(87,121){\footnotesize{$\cdots$}}
\put(107,121){\footnotesize{$1$}}
\put(124,121){\footnotesize{$\cdots$}}

\psline{->}(70,10)(90,10)
\put(79,6){$\rm{T}_1$}

\put(11,38){\footnotesize{$(0,1)$}}
\put(11,58){\footnotesize{$(0,2)$}}
\put(14,78){\footnotesize{$\vdots$}}
\put(11,98){\footnotesize{$(0,n)$}}
\put(4,118){\footnotesize{$(0,n+1)$}}

\psline{->}(10, 70)(10,90)
\put(5,80){$\rm{T}_2$}

\put(20,28){\footnotesize{$y$}}
\put(20,48){\footnotesize{$1$}}
\put(20,68){\footnotesize{$1$}}
\put(22,88){\footnotesize{$\vdots$}}
\put(20,108){\footnotesize{$1$}}
\put(22,128){\footnotesize{$\vdots$}}

\put(40,28){\footnotesize{$\frac{ay}{x}$}}
\put(40,48){\footnotesize{$1$}}
\put(40,68){\footnotesize{$1$}}
\put(42,88){\footnotesize{$\vdots$}}
\put(40,108){\footnotesize{$1$}}
\put(42,128){\footnotesize{$\vdots$}}

\put(60,28){\footnotesize{$\frac{ay}{x}$}}
\put(60,48){\footnotesize{$1$}}
\put(60,68){\footnotesize{$1$}}
\put(62,88){\footnotesize{$\vdots$}}
\put(60,108){\footnotesize{$1$}}
\put(62,128){\footnotesize{$\vdots$}}

\put(80,28){\footnotesize{$\frac{ay}{x}$}}
\put(80,48){\footnotesize{$1$}}
\put(80,68){\footnotesize{$1$}}
\put(82,88){\footnotesize{$\vdots$}}
\put(80,108){\footnotesize{$1$}}
\put(82,128){\footnotesize{$\vdots$}}

\put(100,28){\footnotesize{$\frac{ay}{x}$}}
\put(100,48){\footnotesize{$1$}}
\put(100,68){\footnotesize{$1$}}
\put(102,88){\footnotesize{$\vdots$}}
\put(100,108){\footnotesize{$1$}}
\put(102,128){\footnotesize{$\vdots$}}

\put(122,28){\footnotesize{$\vdots$}}
\put(122,48){\footnotesize{$\vdots$}}
\put(122,68){\footnotesize{$\vdots$}}
\put(122,88){\footnotesize{$\vdots$}}
\put(122,108){\footnotesize{$\vdots$}}
\put(122,128){\footnotesize{$\vdots$}}

\end{picture}
\end{center}
\caption{Weight diagram of the $2$-variable weighted shift in Propositions %
\ref{prophyp} and \ref{propsub} }
\label{Figure 3}
\end{figure}

\begin{proposition}
\label{prophyp}The $2$-variable weighted shift $\mathbf{T}$ given by Figure %
\ref{Figure 3} is hyponormal if and only if $y\leq \min \{\frac{x}{a},x\sqrt{%
\frac{1-x^{2}}{x^{2}-2a^{2}x^{2}+a^{4}}}\}.$
\end{proposition}

\begin{proof}
By the Six-point Test (Theorem \ref{joint hypo}), to show the joint
hyponormality of $\mathbf{T}$ it is enough to check that 
\begin{equation*}
H:=\left( 
\begin{array}{cc}
1-x^{2} & \frac{a^{2}y}{x}-yx \\ 
\frac{a^{2}y}{x}-yx & 1-y^{2}%
\end{array}%
\right) \geq 0.
\end{equation*}%
Since $x<1$, the positivity of $H$ is equivalent to $\det H\geq 0$, i.e., 
\begin{equation*}
(1-x^{2})(1-y^{2})\geq (\frac{a^{2}y}{x}-yx)^{2},
\end{equation*}%
which in turn is equivalent to $y\leq x\sqrt{\frac{1-x^{2}}{%
x^{2}-2a^{2}x^{2}+a^{4}}}$ (observe that $%
x^{2}-2a^{2}x^{2}+a^{4}=x^{2}(1-x^{2})+(x^{2}-a^{2})^{2}>0).$
\end{proof}

\begin{proposition}
\label{propsub}The $2$-variable weighted shift $\mathbf{T}$ given by Figure %
\ref{Figure 3} is subnormal if and only if $y\leq \sqrt{\frac{1-x^{2}}{%
1-a^{2}}}.$
\end{proposition}

\begin{proof}
From Figure \ref{Figure 3}, it is obvious that $\mathbf{T}_{\mathcal{M}%
}\cong (S_{a}\otimes I,I\otimes U_{+})$ (recall that $S_{a}:=shift(a,1,1,...)
$ and $U_{+}$ is the (unweighted) unilateral shift). \ By Lemma \ref%
{lemtensor}, $\mathbf{T}_{\mathcal{M}}$ is subnormal, with Berger measure $%
\mu _{\mathcal{M}}:=[(1-a^{2})\delta _{0}+a^{2}\delta _{1}]\times \delta _{1}
$. \ By Proposition \ref{backext}, 
\begin{eqnarray*}
\mathbf{T}\text{ is subnormal } &\Leftrightarrow &\beta _{00}^{2}\left\| 
\frac{1}{t}\right\| _{L^{1}(\mu _{\mathcal{M}})}(\mu _{\mathcal{M}%
})_{ext}^{X}\leq \nu  \\
&\Leftrightarrow &y^{2}[(1-a^{2})\delta _{0}+a^{2}\delta _{1}]\leq
(1-x^{2})\delta _{0}+x^{2}\delta _{1} \\
&\Leftrightarrow &y^{2}(1-a^{2})\leq 1-x^{2}\text{ and }ay\leq x \\
&\Leftrightarrow &y\leq \min \{\frac{x}{a},\sqrt{\frac{1-x^{2}}{1-a^{2}}}\}%
\text{.}
\end{eqnarray*}
\end{proof}

We summarize the results in Propositions \ref{prophyp} and \ref{propsub} as
follows.

\begin{theorem}
\label{secondthm}The $2$-variable weighted shift $\mathbf{T}$ given by
Figure \ref{Figure 3} is hyponormal and not subnormal if and only if $x>a$
and $\sqrt{\frac{1-x^{2}}{1-a^{2}}}<y\leq x\sqrt{\frac{(1-x^{2})}{%
x^{2}+a^{4}-2a^{2}x^{2}}}$ (see Figure \ref{Figure 4}).
\end{theorem}

\setlength{\unitlength}{1mm} \psset{unit=1mm}

\begin{figure}[th]
\begin{center}
\begin{picture}(120,82)

\pscustom[linewidth=0pt,fillstyle=solid,fillcolor=gray]{
\pscurve(42,50.4)(47,48)(53,10)}
\pscustom[linewidth=0pt,fillstyle=solid,fillcolor=white]{
\pscurve(42,50.4)(48,38)(53,10)}

\psline{->}(10,10)(10,77)
\psline{->}(10,10)(62,75)
\psline[linewidth=2pt](10,10)(42,50)

\pscurve[linestyle=dashed,dash=3pt 2pt](10,71)(20,68)(42,50.4)(48,38)(53,10)
\pscurve[linestyle=dashed,dash=3pt 2pt](10,10)(28,40)(42,50.4)(47,48)(53,10)

\psline[linewidth=2pt]{<-}(46.9,44)(60,50)
\psline[linewidth=2pt]{<-}(40,34)(60,40)
\put(61,49){$\mathbf{T}$ is hyponormal but not subnormal in this region}
\put(61,39){$\mathbf{T}$ is subnormal in this region}
\psline[linewidth=2pt]{<-}(25,64.5)(27,74.5)
\put(24,79){$y={\sqrt{\frac{1-x^2}{1-a^2}}}$}
\psline[linewidth=2pt]{<-}(53.3,24)(63.3,27)
\put(65,27){$y=x{\sqrt{\frac{1-x^2}{x^2+a^4-2a^2x^2}}}$}
\put(6,6){\footnotesize{$(0,0)$}}
\put(40,5){\footnotesize{$(a,0)$}}
\put(42.5,10){\pscircle*(0,0){1}}
\put(1,41){\footnotesize{$(0,a)$}}
\put(10,42.5){\pscircle*(0,0){1}}
\put(52,10){\pscircle*(1,0){1}}
\put(51,5){\footnotesize{$(1,0)$}}
\put(63,75){$y={\frac{x}{a}}$}
\put(7,79){$y$}

\psline{->}(10,10)(75,10)
\psline[linewidth=2pt](10,10)(53,10)
\put(1,50){\footnotesize{$(0,1)$}}
\put(10,50.4){\pscircle*(0,0){1}}
\put(-6,71){\footnotesize{$(0,\sqrt{\frac{1}{1-a^2}})$}}

\put(77,7){$x$}

\end{picture}
\end{center}
\caption{Regions of hyponormality and subnormality for the $2$-variable
weighted shift in Theorem \ref{secondthm}}
\label{Figure 4}
\end{figure}
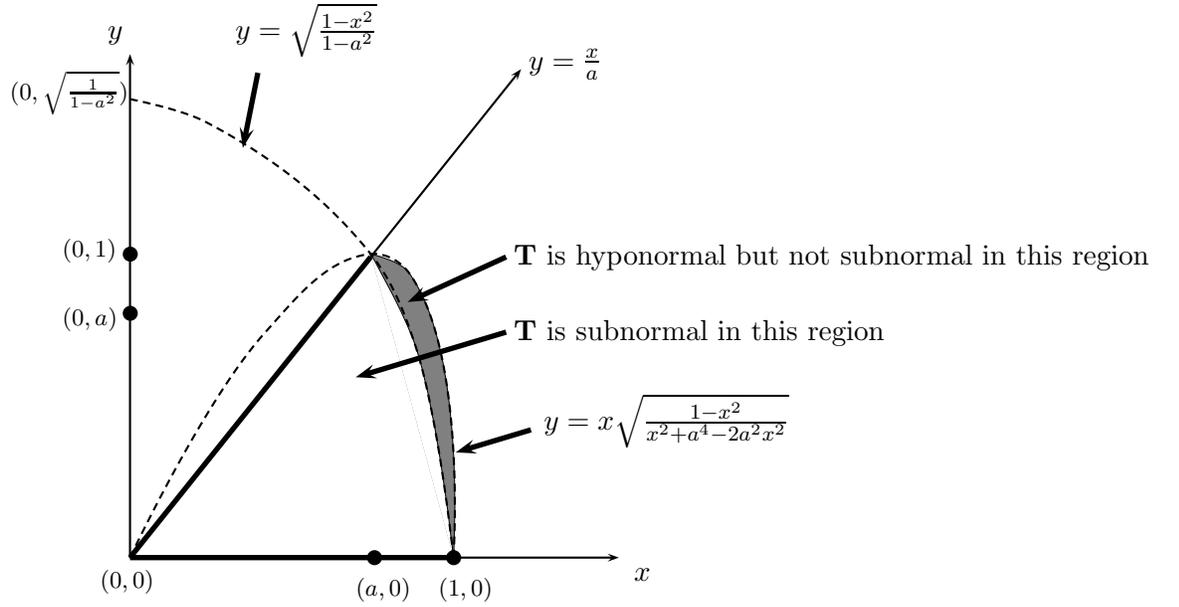

\begin{remark}
As exemplified in Figure \ref{Figure 4}, observe that for $x>a$, $\sqrt{%
\frac{1-x^{2}}{1-a^{2}}}<x\sqrt{\frac{1-x^{2}}{x^{2}+a^{4}-2a^{2}x^{2}}}<%
\frac{x}{a}$; for, if $a<x$ we have 
\begin{eqnarray*}
a^{4} &<&a^{2}x^{2}\Rightarrow x^{2}+a^{4}-2a^{2}x^{2}<(1-a^{2})x^{2} \\
&\Rightarrow &\frac{1-x^{2}}{1-a^{2}}<\frac{x^{2}(1-x^{2})}{%
x^{2}+a^{4}-2a^{2}x^{2}}\text{,}
\end{eqnarray*}%
and 
\begin{eqnarray*}
a^{2}(1-a^{2}) &<&x^{2}(1-a^{2})\Rightarrow a^{2}+a^{2}x^{2}<x^{2}+a^{4} \\
&\Rightarrow &a^{2}(1-x^{2})<x^{2}+a^{4}-2a^{2}x^{2} \\
&\Rightarrow &\frac{x^{2}(1-x^{2})}{x^{2}+a^{4}-2a^{2}x^{2}}<\frac{x^{2}}{%
a^{2}},
\end{eqnarray*}%
as desired.
\end{remark}

\section{\label{fourth}The Third Family of Counterexamples}

\textbf{Construction of the family}. \ Let us consider the following $2$%
-variable weighted shift (see Figure \ref{Figure 5}), where 
\begin{equation}
\left\{ 
\begin{array}{l}
\text{(i) \ }0<\xi _{1}<\xi _{2}<...<\xi _{n}\nearrow 1\text{;} \\ 
\\ 
\text{(ii) \ }W_{\xi }:=shift(\xi _{1},\xi _{2},...)\text{ is subnormal with
Berger measure }\nu \text{;} \\ 
\\ 
\text{(iii) \ }\frac{1}{s^{2}}\in L^{1}(\nu )\text{ (this implies that }%
\frac{1}{s}\in L^{1}(\nu )\text{, by Jensen's inequality);} \\ 
\\ 
\text{(iv) \ }\xi _{e}\equiv \xi _{ext}:=(\int \frac{1}{s}d\nu (s))^{-1/2}%
\text{;} \\ 
\\ 
\text{(v) \ }a\leq \frac{1}{\xi _{e}}(\int \frac{1}{s^{2}}d\nu (s))^{-1/2}%
\text{;} \\ 
\\ 
\text{(vi) \ }b\leq \xi _{e}^{2}\text{ (this implies the condition }b<\xi
_{e}\text{); and} \\ 
\\ 
\text{(vii) \ }a^{2}\leq \frac{b^{2}+\xi _{e}^{2}}{2}\text{.}%
\end{array}%
\right.  \label{conditions}
\end{equation}%
(Recall that $\xi _{e}$ is the maximum possible value for $\xi _{0}$ in
Proposition \ref{backward}.)

\setlength{\unitlength}{1mm} \psset{unit=1mm}

\begin{figure}[th]
\begin{center}
\begin{picture}(140,138)

\psline{->}(20,20)(135,20)
\psline(20,40)(125,40)
\psline(20,60)(125,60)
\psline(20,80)(125,80)
\psline(20,100)(125,100)
\psline(20,120)(125,120)
\psline{->}(20,20)(20,135)
\psline(40,20)(40,125)
\psline(60,20)(60,125)
\psline(80,20)(80,125)
\psline(100,20)(100,125)
\psline(120,20)(120,125)

\put(11,16){\footnotesize{$(0,0)$}}
\put(35,16){\footnotesize{$(1,0)$}}
\put(55,16){\footnotesize{$(2,0)$}}
\put(78,16){\footnotesize{$\cdots$}}
\put(95,16){\footnotesize{$(n,0)$}}
\put(115,16){\footnotesize{$(n+1,0)$}}

\put(27,21){\footnotesize{$a$}}
\put(47,21){\footnotesize{$\xi_{e}$}}
\put(67,21){\footnotesize{$\xi_{1}$}}
\put(87,21){\footnotesize{$\cdots$}}
\put(107,21){\footnotesize{$\xi_{n-1}$}}
\put(124,21){\footnotesize{$\cdots$}}

\put(27,41){\footnotesize{$b$}}
\put(47,41){\footnotesize{$\xi_{1}$}}
\put(67,41){\footnotesize{$\xi_{2}$}}
\put(87,41){\footnotesize{$\cdots$}}
\put(107,41){\footnotesize{$\xi_{n}$}}
\put(124,41){\footnotesize{$\cdots$}}

\put(27,62){\footnotesize{$\frac{\xi_{1}b}{\xi_{e}}$}}
\put(47,61){\footnotesize{$\xi_{1}$}}
\put(67,61){\footnotesize{$\xi_{2}$}}
\put(87,61){\footnotesize{$\cdots$}}
\put(107,61){\footnotesize{$\xi_{n}$}}
\put(124,61){\footnotesize{$\cdots$}}

\put(27,82){\footnotesize{$\frac{\xi_{2}b}{\xi_{e}}$}}
\put(47,81){\footnotesize{$\xi_{1}$}}
\put(67,81){\footnotesize{$\xi_{2}$}}
\put(87,81){\footnotesize{$\cdots$}}
\put(107,81){\footnotesize{$\xi_{n}$}}
\put(124,81){\footnotesize{$\cdots$}}

\put(27,101){\footnotesize{$\cdots$}}
\put(47,101){\footnotesize{$\cdots$}}
\put(67,101){\footnotesize{$\cdots$}}
\put(87,101){\footnotesize{$\cdots$}}
\put(107,101){\footnotesize{$\cdots$}}
\put(124,101){\footnotesize{$\cdots$}}

\put(27,122){\footnotesize{$\frac{\xi_{n-1}b}{\xi_{e}}$}}
\put(47,121){\footnotesize{$\xi_{1}$}}
\put(67,121){\footnotesize{$\xi_{2}$}}
\put(87,121){\footnotesize{$\cdots$}}
\put(107,121){\footnotesize{$\xi_{n}$}}
\put(124,121){\footnotesize{$\cdots$}}

\psline{->}(70,10)(90,10)
\put(79,6){$\rm{T}_1$}

\put(11,38){\footnotesize{$(0,1)$}}
\put(11,58){\footnotesize{$(0,2)$}}
\put(14,78){\footnotesize{$\vdots$}}
\put(11,98){\footnotesize{$(0,n)$}}
\put(4,118){\footnotesize{$(0,n+1)$}}

\psline{->}(10, 70)(10,90)
\put(5,80){$\rm{T}_2$}

\put(20,28){\footnotesize{$a$}}
\put(20,48){\footnotesize{$\xi_{e}$}}
\put(20,68){\footnotesize{$\xi_{1}$}}
\put(22,88){\footnotesize{$\vdots$}}
\put(20,108){\footnotesize{$\xi_{n-1}$}}
\put(22,128){\footnotesize{$\vdots$}}

\put(40,28){\footnotesize{$b$}}
\put(40,48){\footnotesize{$\xi_{1}$}}
\put(40,68){\footnotesize{$\xi_{2}$}}
\put(42,88){\footnotesize{$\vdots$}}
\put(40,108){\footnotesize{$\xi_{n}$}}
\put(42,128){\footnotesize{$\vdots$}}

\put(60,28){\footnotesize{$\frac{\xi_{1}b}{\xi_{e}}$}}
\put(60,48){\footnotesize{$\xi_{1}$}}
\put(60,68){\footnotesize{$\xi_{2}$}}
\put(62,88){\footnotesize{$\vdots$}}
\put(60,108){\footnotesize{$\xi_{n}$}}
\put(62,128){\footnotesize{$\vdots$}}

\put(80,28){\footnotesize{$\frac{\xi_{2}b}{\xi_{e}}$}}
\put(80,48){\footnotesize{$\xi_{1}$}}
\put(80,68){\footnotesize{$\xi_{2}$}}
\put(82,88){\footnotesize{$\vdots$}}
\put(80,108){\footnotesize{$\xi_{n}$}}
\put(82,128){\footnotesize{$\vdots$}}

\put(100,28){\footnotesize{$\frac{\xi_{n-1}b}{\xi_{e}}$}}
\put(100,48){\footnotesize{$\xi_{1}$}}
\put(100,68){\footnotesize{$\xi_{2}$}}
\put(102,88){\footnotesize{$\vdots$}}
\put(100,108){\footnotesize{$\xi_{n}$}}
\put(102,128){\footnotesize{$\vdots$}}

\put(122,28){\footnotesize{$\vdots$}}
\put(122,48){\footnotesize{$\vdots$}}
\put(122,68){\footnotesize{$\vdots$}}
\put(122,88){\footnotesize{$\vdots$}}
\put(122,108){\footnotesize{$\vdots$}}
\put(122,128){\footnotesize{$\vdots$}}

\end{picture}
\end{center}
\caption{Weight diagram of the 2-variable weighted shift in Example 19}
\label{Figure 5}
\end{figure}
Observe that $T_{1}\cong T_{2}$ and that $T_{1}T_{2}=T_{2}T_{1}$. \ We claim
that $T_{1}$ (and therefore $T_{2})$ is subnormal. \ For, the choice of $\xi
_{e}$ immediately implies that $shift(\xi _{e},\xi _{1},\xi _{2},...)$ is
subnormal, with Berger measure $d\nu _{e}(s):=\frac{\xi _{e}^{2}}{s}d\nu (s)$
(cf. Proposition \ref{backward}). \ Another application of Proposition \ref%
{backward} shows that $shift(a,\xi _{e},\xi _{1},...)$ is subnormal if and
only if $\frac{1}{s}\in L^{1}(\nu _{e})$ (i.e., $\frac{1}{s^{2}}\in
L^{1}(\nu )$, which is true by (\ref{conditions})(iii)) and $a^{2}\xi
_{e}^{2}\int \frac{1}{s^{2}}d\nu (s)\leq 1$, which holds by (\ref{conditions}%
)(v)). \ This implies that the restriction of $T_{1}$ to $\bigvee
\{e_{(i,0)}:i\geq 0\}$ is subnormal. \ Moreover, the subnormality of $T_{1}$
when restricted to $\bigvee \{e_{(i,j)}:i\geq 0\}\;(j>0)$ requires that $%
b\leq \xi _{e}$, which holds by (\ref{conditions})(vi)$.$

For a concrete numerical example, consider the probability measure $d\nu
(s):=3s^{2}ds$ on the interval $[0,1]$. \ The measure $\nu $ corresponds to
a subnormal weighted shift with weights $\xi _{1}=\sqrt{\frac{3}{4}}$, $\xi
_{2}=\sqrt{\frac{4}{5}}$, $\xi _{3}=\sqrt{\frac{5}{6}}$, $...$ ; indeed, in
this case $W_{\xi }$ is the restriction of the Bergman shift $B_{+}$ to the
invariant subspace $\mathcal{M}_{2}$ obtained by removing the first two
basis vectors in the canonical orthonormal basis of $\ell ^{2}(\mathbb{Z}%
_{+})$. \ Clearly $\frac{1}{s^{2}}\in L^{1}(\nu )$, and $\int \frac{1}{s^{2}}%
d\nu (s)=3$; moreover, $\int \frac{1}{s}d\nu (s)=\frac{3}{2}$, so in this
case $\xi _{e}=\sqrt{\frac{2}{3}}$. \ Choosing $a=\sqrt{\frac{1}{2}}$ and $b=%
\sqrt{\frac{1}{3}}$ we see that all conditions in (\ref{conditions}) are
satisfied (cf. Corollary \ref{example3}).

\begin{proposition}
\label{hyponormal3}The $2$-variable weighted shift $\mathbf{T}$ given by
Figure \ref{Figure 5} is hyponormal.
\end{proposition}

\begin{proof}
Since the restriction of $\mathbf{T}$ to $\bigvee \{e_{(i,j)}:i,j\geq 1\}$
is clearly subnormal (being unitarily equivalent to $W_{\xi }\bigotimes
I,I\bigotimes W_{\xi }$), and since the weight diagram of $\mathbf{T}$ is
symmetric with respect to the diagonal $j=i$, it suffices to apply the
Six-point Test (Theorem \ref{joint hypo}) to $\mathbf{k}=(i,0)$, with $i\geq
0$. $\ $

\textbf{Case 1}: $\ \mathbf{k}=(0,0)$. \ Here we have 
\begin{eqnarray*}
\left( 
\begin{array}{cc}
\xi _{e}^{2}-a^{2} & b^{2}-a^{2} \\ 
b^{2}-a^{2} & \xi _{e}^{2}-a^{2}%
\end{array}%
\right) &\geq &0\Leftrightarrow (\xi _{e}^{2}-a^{2})^{2}\geq
(b^{2}-a^{2})^{2} \\
&\Leftrightarrow &\xi _{e}^{2}-a^{2}\geq \left| b^{2}-a^{2}\right| .
\end{eqnarray*}%
When $b\leq a$, the last condition is equivalent to $2a^{2}\leq b^{2}+\xi
_{e}^{2}$, which holds by (\ref{conditions})(vii). \ When $b>a$, the
condition is equivalent to $\xi _{e}\geq b$, which is guaranteed by (\ref%
{conditions})(vi).

\textbf{Case 2}: \ $\mathbf{k}=(1,0)$. \ Here 
\begin{eqnarray*}
\left( 
\begin{array}{cc}
\xi _{1}^{2}-\xi _{e}^{2} & \frac{\xi _{1}^{2}b}{\xi _{e}}-b\xi _{e} \\ 
\frac{\xi _{1}^{2}b}{\xi _{e}}-b\xi _{e} & \xi _{1}^{2}-b^{2}%
\end{array}%
\right) &\geq &0\Leftrightarrow (\xi _{1}^{2}-\xi _{e}^{2})(\xi
_{1}^{2}-b^{2})\geq (\frac{\xi _{1}^{2}b}{\xi _{e}}-b\xi _{e})^{2} \\
&\Leftrightarrow &\xi _{1}^{2}-b^{2}\geq (\xi _{1}^{2}-\xi _{e}^{2})\frac{%
b^{2}}{\xi _{e}^{2}}\Leftrightarrow b\leq \xi _{e},
\end{eqnarray*}%
which again is guaranteed by (\ref{conditions})(vi).

\textbf{Case 3}: \ $k=(n+1,0)\;(n\geq 1)$. \ Here 
\begin{eqnarray}
\left( 
\begin{array}{cc}
\xi _{n+1}^{2}-\xi _{n}^{2} & \frac{\xi _{n+1}^{2}b}{\xi _{e}}-\frac{\xi
_{n}^{2}b}{\xi _{e}} \\ 
\frac{\xi _{n+1}^{2}b}{\xi _{e}}-\frac{\xi _{n}^{2}b}{\xi _{e}} & \xi
_{1}^{2}-\frac{\xi _{n}^{2}b^{2}}{\xi _{e}^{2}}%
\end{array}%
\right) &\geq &0  \notag \\
&\Leftrightarrow &(\xi _{n+1}^{2}-\xi _{n}^{2})(\xi _{1}^{2}-\frac{\xi
_{n}^{2}b^{2}}{\xi _{e}^{2}})\geq (\frac{\xi _{n+1}^{2}b}{\xi _{e}}-\frac{%
\xi _{n}^{2}b}{\xi _{e}})^{2}  \notag \\
&\Leftrightarrow &\frac{(\xi _{n+1}^{2}-\xi _{n}^{2})\left( \xi _{1}^{2}\xi
_{e}^{2}-\xi _{n+1}^{2}b^{2}\right) }{\xi _{e}^{2}}\geq 0  \notag \\
&\Leftrightarrow &b\leq \frac{\xi _{1}\xi _{e}}{\xi _{n+1}}\;(\text{for all }%
n\geq 1).  \label{eq in b}
\end{eqnarray}%
Since the sequence $\{\xi _{n}\}$ increases to $1$, the last inequality in (%
\ref{eq in b}) is equivalent to $b\leq \xi _{1}\xi _{e}$, which holds by (%
\ref{conditions})(vi).

The proof is now complete.
\end{proof}

\begin{proposition}
\label{subnormal3}The $2$-variable weighted shift $\mathbf{T}$ given by
Figure \ref{Figure 5} is not subnormal if $p<0$, where $p:=\xi _{e}^{2}\xi
_{1}^{4}+4a^{2}b^{2}\xi _{1}^{2}-b^{2}\xi _{1}^{4}-a^{2}b^{2}\xi
_{e}^{2}-a^{2}b^{4}-2a^{2}\xi _{1}^{4}$.
\end{proposition}

\begin{proof}
Assume that $\mathbf{T}$ is subnormal, and consider the moment matrix
associated \ to the monomials $1$, $x$, $y$ and $yx$ (cf. \cite{tcmp1}, \cite%
{tcmp6}), that is, 
\begin{equation*}
M:=\left( 
\begin{array}{cccc}
1 & a^{2} & a^{2} & a^{2}b^{2} \\ 
a^{2} & a^{2}\xi _{e}^{2} & a^{2}b^{2} & a^{2}b^{2}\xi _{1}^{2} \\ 
a^{2} & a^{2}b^{2} & a^{2}\xi _{e}^{2} & a^{2}b^{2}\xi _{1}^{2} \\ 
a^{2}b^{2} & a^{2}b^{2}\xi _{1}^{2} & a^{2}b^{2}\xi _{1}^{2} & a^{2}b^{2}\xi
_{1}^{4}%
\end{array}%
\right) .
\end{equation*}%
In the presence of a representing measure, it is well known that $M$ must be
positive semi-definite, so in particular $\det M\geq 0$. \ Now, a
straightforward calculation shows that 
\begin{equation*}
\det M=a^{6}b^{2}\left( \xi _{e}^{2}-b^{2}\right) \left( \xi _{e}^{2}\xi
_{1}^{4}-\xi _{e}^{2}a^{2}b^{2}-2a^{2}\xi _{1}^{4}-b^{2}\xi
_{1}^{4}+4a^{2}b^{2}\xi _{1}^{2}-b^{4}a^{2}\right) =a^{6}b^{2}\left( \xi
_{e}^{2}-b^{2}\right) p.
\end{equation*}%
It follows that $p\geq 0$. \ Therefore, $\mathbf{T}$ is not subnormal
whenever $p<0$, as desired.
\end{proof}

\begin{theorem}
\label{theorem3}Let $a>0$ be such that $\sqrt{\frac{\xi _{e}^{2}}{2}}<a\leq 
\sqrt{\frac{\xi _{e}^{2}+\xi _{e}^{4}}{2}}$ and $a\leq \frac{1}{\xi _{e}}%
(\int \frac{1}{s^{2}}d\nu (s))^{-1/2}$, and define $b:=\sqrt{2a^{2}-\xi
_{e}^{2}}$. \ Then the $2$-variable weighted shift $\mathbf{T}\equiv \mathbf{%
(}T_{1},T_{2})$ satisfies (\ref{conditions})(i)-(vii), is hyponormal, and is
not subnormal.
\end{theorem}

\begin{proof}
Observe that the condition $\sqrt{\frac{\xi _{e}^{2}}{2}}<a$ guarantees that 
$2a^{2}>\xi _{e}^{2}$ (so $b$ is well defined), and that the condition $%
a\leq \sqrt{\frac{\xi _{e}^{2}+\xi _{e}^{4}}{2}}$ is equivalent to $%
2a^{2}-\xi _{e}^{2}\leq \xi _{e}^{4}$ (so $b$ satisfies (\ref{conditions}%
)(vi)). \ Moreover, $a^{2}=\frac{b^{2}+\xi _{e}^{2}}{2}$ trivially, so (\ref%
{conditions})(vii) also holds. \ It follows that $T$ is hyponormal, by
Proposition \ref{hyponormal3}. \ To break subnormality, by Proposition \ref%
{subnormal3} it suffices to show that $p$ is negative. \ Since $%
b^{2}=2a^{2}-\xi _{e}^{2}$, we have 
\begin{eqnarray*}
p &=&\xi _{e}^{2}\xi _{1}^{4}-\xi _{e}^{2}a^{2}(2a^{2}-\xi
_{e}^{2})-2a^{2}\xi _{1}^{4}-(2a^{2}-\xi _{e}^{2})\xi
_{1}^{4}+4a^{2}(2a^{2}-\xi _{e}^{2})\xi _{1}^{2}-(2a^{2}-\xi
_{e}^{2})^{2}a^{2} \\
&=&-2\left( \xi _{1}^{2}-a^{2}\right) ^{2}\left( 2a^{2}-\xi _{e}^{2}\right)
<0,
\end{eqnarray*}%
as desired. \ The proof is now complete.
\end{proof}

\begin{corollary}
\label{example3}\cite{DrMcC}Let $dv(s):=3t^{2}ds$ on $[0,1]$ and choose $a=%
\sqrt{\frac{1}{2}}$ and $b=\sqrt{\frac{1}{3}}$. \ Then the $2$-variable
weighted shift $\mathbf{T}$ given by Figure \ref{Figure 5} is commuting, has
subnormal components, is hyponormal, but is not subnormal.
\end{corollary}

\begin{proof}
By Theorem \ref{theorem3} and the comments preceding Proposition \ref%
{hyponormal3}, it suffices to verify that $\sqrt{\frac{\xi _{e}^{2}}{2}}%
<a\leq \sqrt{\frac{\xi _{e}^{2}+\xi _{e}^{4}}{2}}$. \ Since $\xi _{e}=\sqrt{%
\frac{2}{3}}$ and $a=\sqrt{\frac{1}{2}}$, the result follows by a
straightforward calculation.
\end{proof}

\section{An Instance When Hyponormality Suffices}

In this section we will prove that under a suitable condition hyponormality
does imply subnormality for commuting pairs of subnormal operators. \ We
begin with an elementary result of independent interest.

\begin{lemma}
\label{lemspecial}Let $\nu $ be a probability measure on $[0,1]$, and let $%
\gamma _{n}\equiv \gamma _{n}(\nu ):=\int s^{n}d\nu (s)\;(n\geq 0)$ be the
moments of $\nu $. \ The sequence $\{\gamma _{n}\}_{n=0}^{\infty }$ is
bounded below if and only if $\nu $ has an atom at $\{1\}$.
\end{lemma}

\begin{proof}
$(\Leftarrow )$ \ Let $\rho :=\nu (\{1\})>0$ and write $\nu \equiv (1-\rho
)\eta +\rho \delta _{1}$, where $\eta $ is a probability measure on $[0,1]$
with $\eta (\{1\})=0.$ \ It follows that $\gamma _{n}(\nu )\geq \rho \int
s^{n}\;d\delta _{1}(s)=\rho \;($all $n\geq 0)$, so $\{\gamma _{n}\}$ is
bounded below by $\rho $.

$(\Rightarrow )$ \ Suppose $\nu (\{1\})=0$, let $f_{n}(s):=s^{n}\;(0\leq
s\leq 1$, $n\geq 0)$, and consider the sequence of nonnegative functions $%
\{f_{n}\}_{n\geq 0}$. \ Clearly $f_{n}\searrow \chi _{\{1\}}$ pointwise, and 
$\left| f_{n}\right| \leq 1\;($all $n\geq 0)$. \ By the Lebesgue Dominated
Convergence Theorem, $\lim_{n\rightarrow \infty }\gamma
_{n}=\lim_{n\rightarrow \infty }\int s^{n}d\nu (s)=\lim_{n\rightarrow \infty
}s^{n}d\nu (s)=\nu (\{1\})=0$. $\ $Therefore, $\{\gamma _{n}\}$ is not
bounded below.
\end{proof}

We now consider the $2$-variable weighted shift $\mathbf{T}$ given by Figure %
\ref{Figure 6}), where $W_{\xi }:=shift(\xi _{0},\xi _{1},\cdots )$ is a
subnormal contraction with associated measure $\nu $, and $y\leq 1$. \ 
\newline
\setlength{\unitlength}{1mm} \psset{unit=1mm}

\begin{figure}[th]
\begin{center}
\begin{picture}(140,138)

\psline{->}(20,20)(135,20)
\psline(20,40)(125,40)
\psline(20,60)(125,60)
\psline(20,80)(125,80)
\psline(20,100)(125,100)
\psline(20,120)(125,120)
\psline{->}(20,20)(20,135)
\psline(40,20)(40,125)
\psline(60,20)(60,125)
\psline(80,20)(80,125)
\psline(100,20)(100,125)
\psline(120,20)(120,125)

\put(11,16){\footnotesize{$(0,0)$}}
\put(35,16){\footnotesize{$(1,0)$}}
\put(55,16){\footnotesize{$(2,0)$}}
\put(78,16){\footnotesize{$\cdots$}}
\put(95,16){\footnotesize{$(n,0)$}}
\put(115,16){\footnotesize{$(n+1,0)$}}

\put(27,21){\footnotesize{$\xi_{0}$}}
\put(47,21){\footnotesize{$\xi_{1}$}}
\put(67,21){\footnotesize{$\xi_{2}$}}
\put(87,21){\footnotesize{$\cdots$}}
\put(107,21){\footnotesize{$\xi_{n}$}}
\put(124,21){\footnotesize{$\cdots$}}

\put(27,41){\footnotesize{$1$}}
\put(47,41){\footnotesize{$1$}}
\put(67,41){\footnotesize{$1$}}
\put(87,41){\footnotesize{$\cdots$}}
\put(107,41){\footnotesize{$1$}}
\put(124,41){\footnotesize{$\cdots$}}

\put(27,61){\footnotesize{$1$}}
\put(47,61){\footnotesize{$1$}}
\put(67,61){\footnotesize{$1$}}
\put(87,61){\footnotesize{$\cdots$}}
\put(107,61){\footnotesize{$1$}}
\put(124,61){\footnotesize{$\cdots$}}

\put(27,81){\footnotesize{$\cdots$}}
\put(47,81){\footnotesize{$\cdots$}}
\put(67,81){\footnotesize{$\cdots$}}
\put(87,81){\footnotesize{$\cdots$}}
\put(107,81){\footnotesize{$\cdots$}}
\put(124,81){\footnotesize{$\cdots$}}

\put(27,101){\footnotesize{$1$}}
\put(47,101){\footnotesize{$1$}}
\put(67,101){\footnotesize{$1$}}
\put(87,101){\footnotesize{$\cdots$}}
\put(107,101){\footnotesize{$1$}}
\put(124,101){\footnotesize{$\cdots$}}

\put(27,121){\footnotesize{$1$}}
\put(47,121){\footnotesize{$1$}}
\put(67,121){\footnotesize{$1$}}
\put(87,121){\footnotesize{$\cdots$}}
\put(107,121){\footnotesize{$1$}}
\put(124,121){\footnotesize{$\cdots$}}

\psline{->}(70,10)(90,10)
\put(79,6){$\rm{T}_1$}

\put(11,38){\footnotesize{$(0,1)$}}
\put(11,58){\footnotesize{$(0,2)$}}
\put(14,78){\footnotesize{$\vdots$}}
\put(11,98){\footnotesize{$(0,n)$}}
\put(4,118){\footnotesize{$(0,n+1)$}}

\psline{->}(10, 70)(10,90)
\put(5,80){$\rm{T}_2$}

\put(20,28){\footnotesize{$y$}}
\put(20,48){\footnotesize{$1$}}
\put(20,68){\footnotesize{$1$}}
\put(22,88){\footnotesize{$\vdots$}}
\put(20,108){\footnotesize{$1$}}
\put(22,128){\footnotesize{$\vdots$}}

\put(40,28){\footnotesize{$\frac{y}{\xi_{0}}$}}
\put(40,48){\footnotesize{$1$}}
\put(40,68){\footnotesize{$1$}}
\put(42,88){\footnotesize{$\vdots$}}
\put(40,108){\footnotesize{$1$}}
\put(42,128){\footnotesize{$\vdots$}}

\put(60,28){\footnotesize{$\frac{y}{\xi_{0}\xi_{1}}$}}
\put(60,48){\footnotesize{$1$}}
\put(60,68){\footnotesize{$1$}}
\put(62,88){\footnotesize{$\vdots$}}
\put(60,108){\footnotesize{$1$}}
\put(62,128){\footnotesize{$\vdots$}}

\put(100,28){\footnotesize{$\frac{y}{\sqrt{\gamma_n}}$}}
\put(100,48){\footnotesize{$1$}}
\put(100,68){\footnotesize{$1$}}
\put(102,88){\footnotesize{$\vdots$}}
\put(100,108){\footnotesize{$1$}}
\put(102,128){\footnotesize{$\vdots$}}

\put(122,28){\footnotesize{$\vdots$}}
\put(122,48){\footnotesize{$\vdots$}}
\put(122,68){\footnotesize{$\vdots$}}
\put(122,88){\footnotesize{$\vdots$}}
\put(122,108){\footnotesize{$\vdots$}}
\put(122,128){\footnotesize{$\vdots$}}

\end{picture}
\end{center}
\caption{Weight diagram of the 2-variable weighted shift in Theorem \ref%
{hyposubnormal}}
\label{Figure 6}
\end{figure}
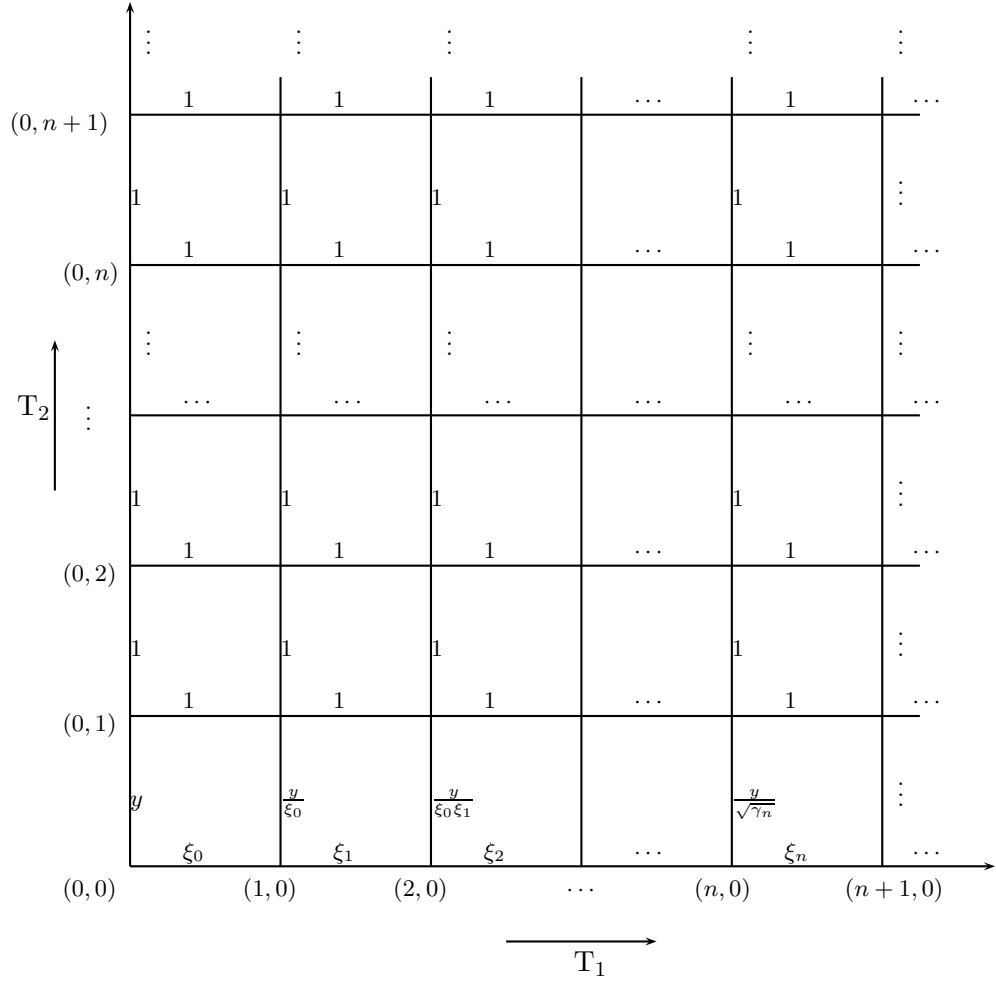
It is clear that $T_{1}T_{2}=T_{2}T_{1}$, and that $T_{1}$ is subnormal
(being the orthogonal direct sum of $W_{\xi }$ and copies of $U_{+}$). \ To
ensure the subnormality of $T_{2}$, we must impose the condition $\frac{y}{%
\sqrt{\gamma _{n}}}\leq 1\;($all $n\geq 0)$, i.e., $y^{2}\leq \gamma _{n}\;($%
all $n\geq 0)$, where $\gamma _{n}\equiv \gamma _{n}(\nu )$. \ Notice that
this condition also guarantees the boundedness of $\mathbf{T}$. \ 

\begin{theorem}
\label{hyposubnormal}Let $\mathbf{T}$ be the $2$-variable weighted shift
given by Figure \ref{Figure 6}, and assume that $\mathbf{T}$ is hyponormal.
\ Then $\mathbf{T}$ is subnormal.
\end{theorem}

\begin{proof}
We apply the Six-point Test (Theorem \ref{joint hypo}) to an arbitrary
lattice point of the form $(n,0)$. \ Since $\mathbf{T}$ is hyponormal by
hypothesis, we must have $(\xi _{n+1}^{2}-\xi _{n}^{2})(1-\frac{y^{2}}{%
\gamma _{n}})\geq (\frac{y}{\sqrt{\gamma _{n+1}}}-\frac{y\xi _{n}}{\sqrt{%
\gamma _{n}}})^{2}$, or equivalently $(\xi _{n+1}^{2}-\xi _{n}^{2})(1-\frac{%
y^{2}}{\gamma _{n}})\geq \frac{y^{2}}{\gamma _{n}}(\frac{1}{\xi _{n}}-\xi
_{n})^{2}$, that is, $y^{2}\leq (\frac{\xi _{n+1}^{2}-\xi _{n}^{2}}{\xi
_{n+1}^{2}+\frac{1}{\xi _{n}^{2}}-2})\gamma _{n}$. \ Since $\xi _{n}^{2}+%
\frac{1}{\xi _{n}^{2}}-2\geq 0$ and $\frac{\xi _{n+1}^{2}-\xi _{n}^{2}}{\xi
_{n+1}^{2}+\frac{1}{\xi _{n}^{2}}-2}=\frac{\xi _{n+1}^{2}-\xi _{n}^{2}}{(\xi
_{n+1}^{2}-\xi _{n}^{2})+\xi _{n}^{2}+\frac{1}{\xi _{n}^{2}}-2}$, it follows
that $\frac{\xi _{n+1}^{2}-\xi _{n}^{2}}{\xi _{n+1}^{2}+\frac{1}{\xi _{n}^{2}%
}-2}\leq 1,$ so $0<y^{2}\leq \gamma _{n}\;($all $n\geq 0)$. \ Thus, $%
\{\gamma _{n}\}$ is bounded below, and by Lemma \ref{lemspecial} we can
write\ $\nu =(1-\rho )\eta +\rho \delta _{1}$, with $\rho :=\nu (\{1\})$ and 
$\eta (\{1\})=0$. \ It follows that $y^{2}\leq \rho $. $\ $Thus, $%
y^{2}\delta _{1}\leq \nu .$ \ By Proposition \ref{backext}, $\mathbf{T}$ is
subnormal.
\end{proof}

\begin{remark}
Theorem \ref{hyposubnormal} (and its proof) reveals that for the $2$%
-variable weighted shift given by Figure \ref{Figure 6}, the subnormality of 
$T_{2}$ is equivalent to the subnormality of $\mathbf{T}$, which in turn is
equivalent to the hyponormality of $\mathbf{T}$.
\end{remark}

\end{document}